\newif\ifArXiV

\ArXiVtrue	%ArXiV style
\ifArXiV
\documentclass[a4paper]{article}
\usepackage[affil-it]{authblk}
\usepackage[authoryear]{natbib}
\newenvironment{frontmatter}{}{}

\let\address\affil
\newenvironment{keyword}{\small \textbf{Keywords: }}{}
\else
\documentclass[authoryear,preprint]{elsarticle}     % onecolumn (second format)
\fi
\usepackage{fullpage}

\usepackage{graphicx}
\usepackage{pdflscape}
\usepackage{secdot}
\usepackage{subcaption}
\usepackage{amsmath, amssymb,amsthm}
\usepackage{amsfonts}
\usepackage{latexsym}
\usepackage{booktabs}
\usepackage{enumitem,blindtext}
\usepackage[export]{adjustbox}
\usepackage{mathrsfs}
\usepackage[usenames,dvipsnames]{color}
\usepackage[mathcal]{eucal}
\usepackage{caption}
\usepackage{multirow}
\usepackage{amsmath,amsfonts}
\usepackage[linesnumbered, noend, ruled]{algorithm2e}
\usepackage{setspace}
\newcommand{\comments}[1]{}
\DeclareGraphicsExtensions{.png,.pdf,.jpg}

\newtheorem{theorem}{Theorem}
\newtheorem{definition}{Definition}

\newtheorem{proposition}{Proposition}

\newtheorem{remark}{Remark}
\newtheorem{example}{Example}

\usepackage{tikz}
\usepackage{tikzsymbols}
\usetikzlibrary{arrows,automata}
\usetikzlibrary{shapes.geometric, arrows}
\usetikzlibrary{positioning}

 \newcommand{\KT}[1]{{\color{black}#1}} 
 \newcommand{\MS}[1]{{\color{black}#1}} 
 \newcommand{\NA }[1]{{\color{black}#1}} 
\newcommand{\EG}[1]{{\color{black}#1}} 

\newcommand{\RMS}[1]{{\color{black}#1}} 
\newcommand{\RKT}[1]{{\color{black}#1}}
\newcommand{\RNA}[1]{{\color{black}#1}}

\newcommand{\DSPa}[0]{$\text{DSP1}(\bar{x},\omega)$\xspace}
\newcommand{\SPa}[0]{$\text{SP1}(\bar{x},\omega)$\xspace}

\newcommand{\Problem}[0]{MBSMP\xspace}
\newcommand{\BC}[0]{BC\xspace}
\newcommand{\ILP}[0]{ILP\xspace}
\newcommand{\BD}[0]{BD\xspace}
\newcommand{\BDs}[0]{BDs\xspace}
\newcommand{\IC}[0]{IC\xspace}

\usepackage[normalem]{ulem}
\DeclareMathOperator*{\argmin}{arg\,min}

\usepackage{xstring}

\usepackage{hyperref}

\newcommand{\mytag}[1]{(\hypertarget{#1}{\mathrm{#1}})}
\newcommand{\myref}[1]{\textnormal{(\hyperlink{#1}{#1})}}

\begin{document}

\begin{frontmatter}

\title{Benders decomposition algorithms for minimizing the spread of harmful contagions in networks}

\ifArXiV
\author[1]{K\"ubra Tan{\i}nm{\i}\c{s}\thanks{ktaninmis@ku.edu.tr}}
\author[2]{Necati Aras\thanks{arasn@bogazici.edu.tr}}
\author[3]{Evren G\"uney\thanks{guneye@mef.edu.tr}}
\author[4]{Markus Sinnl\thanks{markus.sinnl@jku.at}}
\date{}
\else
\author[1]{K\"ubra Tan{\i}nm{\i}\c{s}}
\ead{ktaninmis@ku.edu.tr}
\author[2]{Necati Aras}
\ead{arasn@bogazici.edu.tr}
\author[3]{Evren G\"uney}
\ead{guneye@mef.edu.tr}
\author[4]{Markus Sinnl}
\ead{markus.sinnl@jku.at}
\fi

%\author{Other Author 3}
\address[1]{Department of Industrial Engineering, Ko\c{c} University, \.{I}stanbul, Turkey }
\address[2]{Department of Industrial Engineering, Bo\u{g}azi\c{c}i University, \.{I}stanbul, Turkey }
\address[3]{Department of Industrial Engineering, MEF University, Turkey }
\address[4]{Institute of Business Analytics and Technology Transformation/JKU Business School, Johannes Kepler University Linz, Austria}
\ifArXiV
\maketitle	
\fi
\begin{abstract}

The COVID-19 pandemic has been a recent example for the spread of a harmful contagion in large populations. Moreover, the spread of harmful contagions is not only restricted to an infectious disease, but \NA {is also relevant to} computer viruses and malware in computer networks. Furthermore, the spread of fake news and propaganda in online social networks is also of \MS{major} concern. In this study, we introduce the measure-based spread minimization problem (\Problem), which can help policy makers in minimizing the spread of harmful contagions in large networks. We develop exact solution methods based on branch-and-Benders-cut algorithms that make use of the application of Benders decomposition method to two different mixed-integer programming formulations of the \Problem: an arc-based formulation and a path-based formulation. We show that for both formulations the Benders optimality cuts can be generated using a combinatorial procedure rather than solving the dual subproblems using linear programming. Additional improvements such as using scenario-dependent extended seed sets, initial cuts, and a starting heuristic are also incorporated into our branch-and-Benders-cut algorithms. We investigate the contribution of various components of the solution algorithms to the performance on the basis of computational results obtained on a set of instances derived from existing ones in the literature.
\end{abstract}

\ifArXiV

\else

\fi

\end{frontmatter}

\section{Introduction and motivation}

The COVID-19 pandemic is a somber reminder of the danger of the spread of harmful contagions in networks. It has infected over 759 million people and caused the death of almost 6.8 million people\footnote{according to \url{https://covid19.who.int/}, accessed on March 18, 2023}. The pandemic also caused serious economic damage, see, e.g., \cite{magistretti2021after}. Moreover, the spread of harmful contagions is not only restricted to an infectious disease, but can also occur as computer viruses and malware in computer networks \citep{zimba2019economic}. Furthermore, the spread of fake news and propaganda in online social networks is also a very pressing issue \citep{lazer2018science}. 

In this work, we introduce the \emph{measure-based spread minimization problem} (\Problem), which can be used to model the way how to optimally minimize the spread of harmful contagions in networks. Suppose that we are given a directed graph $G=(V,A)$ representing a network, a stochastic diffusion model for the spread in the network, and a set of initially infected nodes $I \subset V$ from which the spread of the contagion starts in the network according to the diffusion model. We also have a set of labels $K$ for the arcs, each of which represents a certain relationship (contact) type. While there can be parallel arcs between a node pair, each arc is categorized with exactly one label. Blocking a label means taking a measure that prevents the corresponding contact between every pair of nodes connected via an arc having that label. In other words, there is a measure associated with each label, and taking a measure implies removing the arcs with the corresponding label from the network. In the context of disease spread, some possible measures could be closing of schools, closing of department stores, closing of restaurants, and the lock-down of a certain area. The goal is to take measures (i.e., remove sets of arcs) in such a way that the spread of the harmful contagion is minimized while respecting a given budget for taking the measures. In this work, we focus on the \emph{independent cascade model (\IC)} as the stochastic diffusion model. This is a popular diffusion model which has been used for many different optimization problems related to spread or influence maximization/minimization in the context of (social) networks. In the \IC model, each infected node has a ``one-shot'' chance to spread the contagion. We discuss the \IC model in more detail in Section \ref{sec:deftheory}, where we also give a formal definition of the \Problem.

The \Problem is related to the \emph{spread blocking problems} considered in \cite{kuhlman2013blocking} \NA{among others}. In \cite{kuhlman2013blocking}, the authors consider a deterministic diffusion model and define various problems based on different objective functions. Next to minimizing the spread, they also \RNA{study} the objective of saving all non-seed nodes from infection while minimizing the cost of the needed blocking actions. \RNA{In terms of decisions, they deal with blocking of edges as well as blocking of labels or equivalently actions}. %Moreover, some of the problems consider blocking of edges, while others consider \MS{blocking of} labels. 
\RKT{The authors prove the NP-hardness of both spread minimization via edge blockings and the problem of saving all salvageable nodes via action  blockings, under a complex contagion model. They also show that the problem of saving all salvageable nodes via edge blockings is efficiently solvable.}
%The authors prove the NP-hardness for some of the resulting problems, while for other variants they show that they can be solved in polynomial time. 
For the two NP-hard problems, \cite{kuhlman2013blocking} present simple heuristics. 
%Similar heuristics for other variants of spread-blocking problems are also developed for other problem variants
\RKT{Similar heuristics are also developed for several variants of the spread-blocking problems} 
based on edge/node deletion in deterministic networks \citep{eubank2006structure,kuhlman2010finding, enns2015link, kimura2009blocking}. In \citet{gillen2018critical} a deterministic edge-based spread-blocking problem is studied and integer linear programming (\ILP) approaches to solve the problem to proven optimality are presented.  A more detailed review on previous related work can be found in Section \ref{sec:prev}.

\paragraph{Contribution and outline}

\RMS{The contribution of our work consists of the following:
\begin{itemize}
\item We introduce the \Problem, which extends the existing problem in the area of spread minimization by considering i) a stochastic diffusion model, and ii) label-based blocking of arcs. \RMS{This is a combination which, to the best of our knowledge, has not been considered before in literature.}
\item  We show that contrary to many other spread-related problems such as the influence maximization problem, the objective function in our problem is neither \emph{supermodular} nor \emph{submodular}.
\item We present two integer programming formulations for the \Problem, which allow us to solve the problem to proven optimality, while most of the work in the area focuses on heuristic approaches.
\item We show how to apply Benders decomposition (\BD) to both of our formulations. \RNA{This gives rise} to obtain models with less decision variables, which allows for better scalability of our solution algorithms.
\item We show that the Benders optimality cuts in both \BDs can be separated using a combinatorial procedure, i.e., there is no need to solve linear programming subproblems to obtain the optimality cuts.
\item We discuss how the the second formulation can be strengthened by down-lifting inequalities present in the \ILP and how this strengthening can be incorporated in the \BD. 
\item We design Branch-and-Benders-cut (\BC) solution algorithms, where the Benders optimality cuts are generated as needed. We also develop additional improvements to the solution algorithms, which include using i) scenario-dependent extended seed sets, ii) initial cuts, and iii) a starting heuristic.
\item We present a computational study on synthetic instances and \RNA{benchmark} network instances from the literature. In the computational study, we analyze the influence of both the various components of our solution algorithm and instance characteristics on the solution. 
\end{itemize}
}

The remainder of the paper is organized as follows: In Section \ref{sec:prev} we provide a discussion of previous and related work. In Section \ref{sec:deftheory} we give a formal problem definition, including a discussion of the \IC used to model the diffusion. We also show that the objective function in our problem is neither \emph{supermodular} nor \emph{submodular}. Section \ref{sec:arc-based} contains the compact arc-based formulation and the associated \BD method and Section \ref{sec:path-based} contains the path-based formulation and the associated \BD method. In this section we also describe how to down-lift the constraints of the path-based formulation, discuss the relationship to the compact formulation from the previous section, and show that the Benders optimality cuts can be separated in a combinatorial fashion. The details of our \BC algorithms are explained in Section \ref{sec:algorithmic_details}, including a description of our proposed enhancements to the algorithms and implementation details of the cut separation. Computational results are provided in Section \ref{sec:results}. We conclude the paper in Section \ref{sec:conclusion} with some remarks and suggestions for future research.

\section{Previous and related work} \label{sec:prev}

%\KT{TODO: @all: check the recent papers}

\RMS{In the literature there exist several works which study the spread of a harmful contagion in a network. As already discussed in the introduction, there are several aspects which need to be taken into account when dealing with such a problem. One of them is how to model the spread (i.e., in a deterministic or stochastic fashion) and another is how to model the measures that need to be taken against the spread. Moreover, there exist heuristic methods as well as exact solution approaches developed for problems related to spread in networks. We start our literature review by discussing heuristics for spread blocking considering deterministic spread models in Section \ref{sec:prevheur}, followed by exact solution approaches for the same problem setting in Section \ref{sec:prevexact}. To the best of our knowledge, the only previous work which used a stochastic diffusion model in the context of spread minimization is \cite{taninmics2022improved}, where a bilevel stochastic spread-blocking problem based on node-deletion is considered and solved using \ILP. Thus, to provide some context for stochastic diffusion models, in Section \ref{sec:previmp} we give an overview on \emph{influence maximization problems (IMPs)}, where such diffusion models are often used. Another recently emerging type of problems which is concerned with blocking nodes/edges in (social) networks are social network interdiction problems, which are discussed in Section \ref{sec:prevsocial}. However, while these problems make sense in the area of social network analysis, their applicability in preventing the spread of harmful contagions in networks is rather limited, as the resulting networks after interdiction can still be highly connected.

\subsection{Heuristics for spread blocking considering deterministic spread models} \label{sec:prevheur}
A group of studies such as \cite{eubank2006structure}, \citet{kuhlman2010finding,kuhlman2013blocking}%, \citet{kuhlman2013blocking}, 
, \cite{wang2013negative}, and \cite{ju2021node} focus minimizing disease spread via edge/node deletion. They propose greedy algorithms or heuristic methods while considering a deterministic model for the spread of the contagion. Some of these methods are simply based on centrality measures of the nodes in the network. \cite{enns2015link} survey and compare several heuristics for minimizing disease spread in networks, while \cite{holme2004efficient} proposes a heuristic for vaccination in networks. \cite{kimura2009blocking} consider another heuristic method for blocking links in networks using a ``contamination degree" as the criterion. Interestingly, these heuristic studies often give different recommendations regarding which strategies perform the best, which arises from the fact that heuristics as opposed to optimal solution algorithms are employed, and different underlying assumptions about the network are made. 

\subsection{Exact solution approaches for spread blocking considering deterministic spread models} \label{sec:prevexact}

In \cite{gillen2018critical}, the authors study which arcs to remove to minimize the spread under a deterministic linear threshold diffusion model and present \ILP approaches for their problem. Another deterministic node-deletion problem motivated by the spread of influenza-virus is investigated in \cite{charkhgard2018integer}, but the developed \ILP model which is obtained by linearizing a nonlinear formulation of the problem can only handle small-scale problem instances with up to 200 nodes. In \cite{gillen2021fortification}, a robust version of a node-deletion problem is considered and \ILP models are presented for its solution. In another class of such problems, the connectivity of the network is reduced instead of directly minimizing the disease/influence spread. This is the case in \cite{shen2012exact} where node deletions are used to achieve this goal. A similar problem is also tackled with \ILP approaches and heuristics in \cite{arulselvan2009detecting}. Although reducing connectivity could eventually help to reduce the spread, it is not a true measure of the spread. \cite{nandi2016methods} define two ``spread related" deterministic metrics and develop \ILP formulations to optimize these metrics through link removal. A recent survey about node-deletion problems is given by \cite{lalou2018critical}. 

\subsection{Influence maximization problems}\label{sec:previmp}

In the seminal work of \cite{kempe2003maximizing}, the influence maximization problem is introduced to study the effect of viral marketing. In the IMP, we are given a network $G=(V,A)$ with a stochastic diffusion model and an integer number $k$. It involves finding a seed set $S$ of size $k$ in $V$ to start the stochastic diffusion process so as to maximize the expected number of influenced nodes when the process ends. It is an NP-hard problem under several widely used stochastic diffusion models. The work of \citet{kempe2003maximizing} spurred a flurry of activity in the research community, which resulted in the study of many different variants of the problem. In recent years, exact and heuristic solution approaches have emerged to solve variants of the IMP that occur in the analysis of social networks.

The surveys \citet{sun2011survey}, \cite{li2018influence}, and \cite{banerjee2020survey} provide an overview of the developments in heuristic solution approaches for the IMP. There is a less amount of work on solving this problem exactly. In \cite{wu2018two} an \ILP-based approach which exploits the submodularity of the objective function is proposed. \cite{guney2019optimal} considers a Lagrangian relaxation approach based on an \ILP-formulation for the IMP, while \cite{guney2021large} develop a \BD algorithm to solve the same problem. A robust version of the IMP is introduced and tackled by \ILP approaches in \cite{nannicini2020exact}. In the \emph{least cost-influence maximization problem}, the goal is to find the minimum number of seeds to influence the complete network. Exact \ILP approaches for this problem were designed in \citet{fischetti2018least} and \citet{gunnecc2020branch}. The \emph{weighted target set selection problem} studied in \citet{raghavan2019branch} can be seen as a deterministic variant of the IMP in which the spread is a function of the number of influenced individuals at the previous time step rather than being dependent on the spread probability.

\subsection{Social network interdiction problems} \label{sec:prevsocial}
In the context of social networks, interdiction-type problems that focus on interdicting \emph{cliques} or clique-like structures have been considered recently. For example, \citet{furini2019maximum} are concerned with removing nodes so as to minimize the maximal clique size in the remaining graph, whereas \citet{furini2021branch} study the same problem with edge removals. The weighted versions of such problems are studied in \citet{nasirian2019exact} and \cite{pajouh2019minimum}. }

\section{Problem definition and theoretical results \label{sec:deftheory}}

In this section, we first give an overview on stochastic diffusion models and show how we use the \IC model in the context of our problem. Next, we present a formal problem definition of the \Problem. We then show that contrary to many other influence-spread problems in networks, the objective function of the \Problem is neither submodular nor supermodular.

\subsection{Stochastic diffusion models \label{sec:diff}}
A stochastic diffusion model allows to model the spread of harmful contagions in networks with probabilistic infection between nodes. The diffusion process proceeds in time steps, and at each time step additional nodes may get infected. There exist two stochastic diffusion model types frequently considered in the literature: \NA{\emph{cascade models} and \emph{threshold models}} (see, e.g., \citet{kempe2003maximizing}, \citet{kempe2015maximizing}, \citet{wu2018two}, \citet{han2018influence}, \citet{taninmics2019influence}, \citet{guney2021large}, \citet{kahr2021benders}, and \citet{taninmics2022improved}). In cascade models, each arc $(i,j) \in A$ is assigned a probability $p_{ij}$ which represents the chance that node $i$ is successful in infecting node $j$. Each node $i$ has a single attempt at activating (infecting) node $j$ during the diffusion process, which occurs immediately after node $i$ gets activated itself. In threshold models each node is assigned a (probabilistic) threshold value and each arc $(i,j)$ is assigned a certain influence value. A node $j$ gets infected during the diffusion process if the sum of the influence values on the arcs incoming from the already infected neighboring nodes exceeds the threshold value of node $j$.

In this work, we use a special case of cascade models that has received considerable attention in the literature, namely the \IC model in which each activation probability is constant. Following other works in the literature (see, e.g., \citet{kempe2003maximizing}, \citet{kempe2015maximizing}, \citet{wu2018two}, \citet{guney2021large}, and \citet{kahr2021benders}) we use a discrete set of scenarios $\Omega$ to model the stochastic diffusion process caused by the \IC model. Each scenario $\omega \in \Omega$ can then be represented using a so-called \emph{live-arc} graph $G^\omega=(V,A^\omega)$, where the presence of $(i,j) \in A^\omega$ indicates that node $i$ is successful in activating node $j$ in scenario $\omega$. This result follows from the equivalence between the IC and \emph{triggering set} models as shown by \citet{kempe2003maximizing}. The set $\Omega$ is obtained by Monte-Carlo sampling based on the given probabilities $p_{ij}$ for the arcs $(i,j) \in A$, i.e., for each scenario $\omega \in \Omega$, $(i,j) \in A^\omega$ is determined by a (biased) coin flip according to probability $p_{ij}$. \RKT{Note that, although we focus on the independent cascade model, the methods we propose admit the usage of any diffusion model for which an equivalent triggering set model exists, i.e., the scenarios can be represented via live-arc graphs.}

Let $A_k$ be the set of arcs with label $k \in K$. Given a seed set $I \subset V$ and a set $K'\subset K$ of blocked labels, the quantity  $\Phi^\omega(I,K')$ represents the number of nodes affected by the harmful contagion in scenario $\omega$. Formally, it is defined as
\begin{equation}
  \Phi^\omega(I,K')=\big | \{v \in V: \exists \text{ a path from at least one node }i \in I \text{ to } v \text{ when considering } A^\omega\setminus \cup_{k \in K'} A_k \} \big |, \tag{IC-S-$\omega$} \label{eq:spread}
\end{equation} i.e., the cardinality of the set of nodes which can be reached from the seed set $I$ using the non-blocked live-arcs for the considered scenario. We note that \citet{kempe2003maximizing} showed that the \emph{linear threshold} model which is a special case of threshold models can also be represented using live-arc graphs.

\subsection{Problem definition\label{sec:def}}
A formal definition of the \Problem using the \IC model with a discrete set of scenarios is given in the following. In the remainder of this paper, we consider this variant of the \Problem, and for the ease of readability, we keep the abbreviation \Problem.

\begin{definition}[Measure-based spread minimization problem (\Problem)]
Let $G=(V,A)$ be a directed graph and $K$ be a finite set of labels. For each $k \in K$, we are given a \emph{measure cost} $c_k \geq 0$, and an arc set $A_k$, i.e., arcs having the label $k$, that are disjoint and satisfy $\cup_k A_k =A$. We are also given a set of initially infected nodes $I\subset V$, and probabilities $p_{ij}$ for each $(i,j) \in A$ indicating the probability that node $i$ is successful in activating (influencing) node $j$. Let $\Omega$ be a set of scenarios obtained by using Monte-Carlo sampling based on probabilities $p_{ij}$. Let $\Phi^\omega(I,K')$ denote the number of nodes to which the harmful contagion spreads in scenario $\omega$ for the seed set $I\subset V$ and a set of blocked labels $K'\subset K$ as defined in equation \eqref{eq:spread}. The \emph{measure-based spread minimization problem} consists of finding a set of measures to take (labels to block) within a budget $B$ so that the expected number of nodes where the harmful contagion has spread from the given seed set $I$ is minimized.
%\footnote{\KT{to me the last part sounds like the number of seed nodes}}. 
Mathematically, it is defined as
		\begin{equation}
		\min_{K' \subset K: c(K^\prime)\leq B} f(K')=\sum_{\omega \in \Omega} \frac{1}{|\Omega|} \Phi^\omega(I,K') \tag{P} \label{eq:P}
		\end{equation}
		where $c(K')=\sum_{k \in K'} c_k\leq B$ is the budget constraint.
\end{definition}

\begin{remark}
The graph $G=(V,A)$ can have parallel arcs with the same head and tail nodes with different labels. Suppose that there are $p$ arcs from $i\in V$ to $j \in V$ with labels $k_1,\ldots, k_p$. Thus, node $i$ can spread the contagion to node $j$ via one of the $p$ arcs, i.e., contact types, independently. In order to prevent $i$ from spreading the contagion to $j$ in a scenario, each of the labels $k_1,\ldots, k_p$ of the live-arcs in that scenario needs to be blocked. \RKT{By allowing parallel arcs, different relationships between a pair of nodes can be modeled via distinct arcs which can be treated independently as the IC model requires.}
\label{remark:multiple}
\end{remark}

The \Problem is a generalization of the simple contagion variant of the \emph{small weighted critical edge set problem (SWCES)}, which is an NP-hard problem studied in \cite{eubank2006structure} and \cite{kuhlman2013blocking}. In the SWCES a deterministic diffusion model is considered.  Moreover, there are no edge labels and each edge can only be directly blocked, i.e., each edge has a unique label.

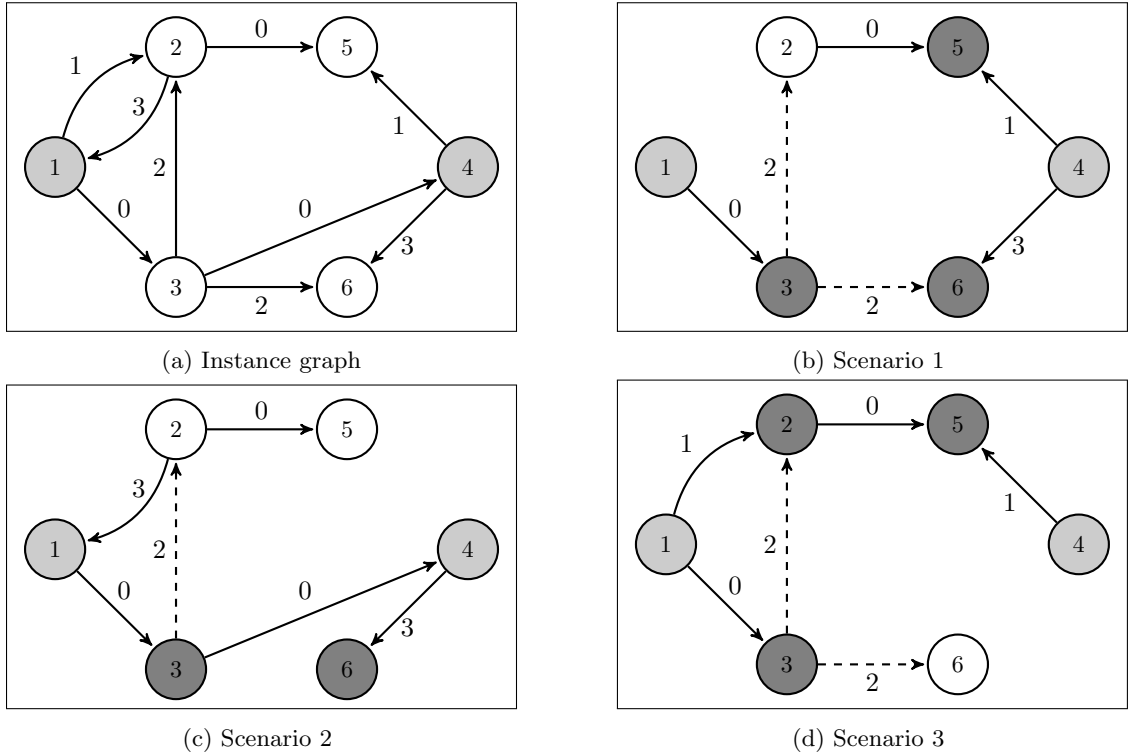
\begin{figure}[h!tb]
\begin{subfigure}{0.5\textwidth}
\centering
 \fbox{
    \begin{tikzpicture}[->,>=stealth',shorten >=1pt,auto,node distance=2.5cm,
                    thick, main node/.style={circle,draw,font=\sffamily\Large\bfseries}]
\tikzstyle{every state}=[fill=white,draw=black,text=black,scale=0.9]

  \node[state,fill=black!20!white] 		   (n1)             		  {1};
  \node[state]         (n2) [above right of=n1] 	  {2};
  \node[state]         (n3) [below right of=n1]       {3};
  \node[state]         (n5) [right of=n2] {5};
  \node[state]         (n0) [right of=n3] {6};
  \node[state,fill=black!20!white]         (n4) [below right of=n5] {4};
  \path (n1) edge[bend left]			  node {1} (n2)
            edge            node[xshift=-3pt] {0} (n3)    
        (n2) edge [bend left]			  node[above, yshift=2pt]  {3} (n1)
        	edge			  node {0} (n5)
        (n3)	edge	   node[left] {2} (n2)
         edge              node {0} (n4)
         edge              node[below] {2} (n0)
        (n4) edge 			  node[xshift=3pt] {1} (n5)
    	    edge			  node[below] {3} (n0)
        ;
\end{tikzpicture}
}
\caption{Instance graph}
\end{subfigure}
\begin{subfigure}{0.5\textwidth}
\centering
 \fbox{
    \begin{tikzpicture}[->,>=stealth',shorten >=1pt,auto,node distance=2.5cm,
                    thick, main node/.style={circle,draw,font=\sffamily\Large\bfseries}]
\tikzstyle{every state}=[fill=white,draw=black,text=black,scale=0.9]

  \node[state,fill=black!20!white] 		   (n1)             		  {1};
  \node[state]         (n2) [above right of=n1] 	  {2};
  \node[state,fill=black!50!white]         (n3) [below right of=n1]       {3};
  \node[state,fill=black!50!white]         (n5) [right of=n2] {5};
  \node[state,fill=black!50!white]         (n0) [right of=n3] {6};
  \node[state,fill=black!20!white]         (n4) [below right of=n5] {4};
  \path (n1)  %edge[bend left]			  node {1} (n2)
        edge            node[xshift=-3pt] {0} (n3)    
        (n2) %edge [bend left]			  node[above, yshift=2pt]  {3} (n1)
        	edge  			  node {0} (n5)
        (n3)	edge [dashed]	   node[left] {2} (n2)
         %edge              node {0} (n4)
         edge     [dashed]         node[below] {2} (n0)
        (n4) edge 			  node[xshift=3pt] {1} (n5)
    	    edge			  node[below] {3} (n0)
        ;
\end{tikzpicture}
}
\caption{Scenario 1 \label{fig:example_s1}}
\end{subfigure}
\vspace*{10pt}
\begin{subfigure}{0.5\textwidth}
\centering
 \fbox{
    \begin{tikzpicture}[->,>=stealth',shorten >=1pt,auto,node distance=2.5cm,
                    thick, main node/.style={circle,draw,font=\sffamily\Large\bfseries}]
\tikzstyle{every state}=[fill=white,draw=black,text=black,scale=0.9]

  \node[state,fill=black!20!white] 		   (n1)             		  {1};
  \node[state]         (n2) [above right of=n1] 	  {2};
  \node[state,fill=black!50!white]         (n3) [below right of=n1]       {3};
  \node[state]         (n5) [right of=n2] {5};
  \node[state,fill=black!50!white]         (n0) [right of=n3] {6};
  \node[state,fill=black!20!white]         (n4) [below right of=n5] {4};
  \path (n1) %edge[bend left]			  node {1} (n2)
            edge            node[xshift=-3pt] {0} (n3)    
        (n2) edge [bend left]			  node[above, yshift=2pt]  {3} (n1)
        	edge			  node {0} (n5)
        (n3)	edge[dashed]	   node[left] {2} (n2)
         edge               node {0} (n4)
        % edge              node[below] {2} (n0)
        (n4) %edge 			  node[xshift=3pt] {1} (n5)
    	    edge			  node[below] {3} (n0)
        ;
\end{tikzpicture}
}
\caption{Scenario 2 \label{fig:example_s2}}
\end{subfigure}
\begin{subfigure}{0.5\textwidth}
\centering
 \fbox{
    \begin{tikzpicture}[->,>=stealth',shorten >=1pt,auto,node distance=2.5cm,
                    thick, main node/.style={circle,draw,font=\sffamily\Large\bfseries}]
\tikzstyle{every state}=[fill=white,draw=black,text=black,scale=0.9]

  \node[state,fill=black!20!white] 		   (n1)             		  {1};
  \node[state,fill=black!50!white]         (n2) [above right of=n1] 	  {2};
  \node[state,fill=black!50!white]         (n3) [below right of=n1]       {3};
  \node[state,fill=black!50!white]         (n5) [right of=n2] {5};
  \node[state]         (n0) [right of=n3] {6};
  \node[state,fill=black!20!white]         (n4) [below right of=n5] {4};
  \path (n1) edge[bend left]			  node {1} (n2)
            edge            node[xshift=-3pt] {0} (n3)    
        (n2)% edge [bend left]			  node[above, yshift=2pt]  {3} (n1)
        	edge 			  node {0} (n5)
        (n3)	edge [dashed]	   node[left] {2} (n2)
         %edge              node {0} (n4)
         edge      [dashed]        node[below] {2} (n0)
        (n4) edge 			  node[xshift=3pt] {1} (n5)
    	    %edge			  node[below] {3} (n0)
        ;
\end{tikzpicture}
}
\caption{Scenario 3 \label{fig:example_s3}}
\end{subfigure}
\caption{An instance of \Problem with six nodes and four arc labels, \MS{seed set} $I=\{1,4\}$ (shown in light grey) and $|\Omega|=3$. For the scenarios, the live-arcs are shown in \ref{fig:example_s1}, \ref{fig:example_s2} and \ref{fig:example_s3}. The budget allows to block at most one label. The solution $K'=\{2\}$ is illustrated, with the blocked arcs shown as dashed arrows. The nodes reached by the spread in each scenario (excluding the seed set) are displayed in dark grey.}\label{fig:problem}
 
\end{figure}

In Figure \ref{fig:problem}, a small instance of the \Problem on a graph with six nodes is given, together with a solution. In this instance $I=\{1,4\}$, $|K|=4$, $c_k=1$, for all $k \in K$, $B=1$ and $|\Omega|=3$. The arc-labels are given as numbers next to the arcs in the figure. In the considered solution, label $2$ is blocked, i.e., $K'=\{2\}$, which is indicated with a dashed line for the arcs with this label in the live-arc graphs for the scenarios. The objective values for the scenarios are $\Phi^1(I,K')=5$ (nodes $1,3,4,5,6$ can be reached), $\Phi^2(I,K')=4$ (nodes $1,3,4,6$ can be reached), $\Phi^3(I,K')=5$ (nodes $1,2,3,4,5$ can be reached), resulting in $f(K')=\frac{1}{3} \cdot 5+\frac{1}{3} \cdot 4+\frac{1}{3} \cdot 5=\frac{14}{3}=4.66$. The optimal solution for this instance is to block label $0$. For this optimal solution, we have $\Phi^1(I,K')=4$ (nodes $1,4,5,6$ can be reached), $\Phi^2(I,K')=3$ (nodes $1,4,6$ can be reached), $\Phi^3(I,K')=4$ (nodes $1,2,4,5$ can be reached) resulting in $f(K')=3.66$

\begin{figure}[h!tb]
    \centering 
    \begin{subfigure}{0.49\textwidth}
        \includegraphics[width=\linewidth]{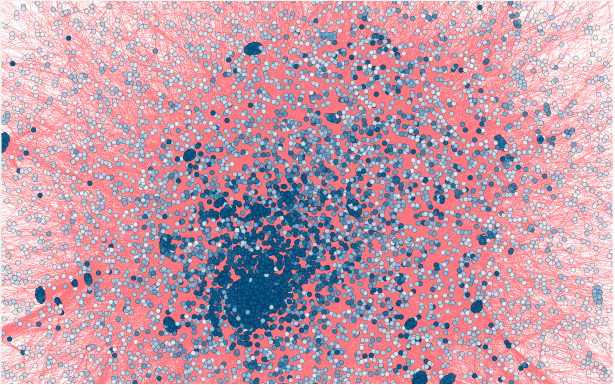}
        \caption{No blocking}
        \label{fig:HepPh_initial}
    \end{subfigure}\hspace{5pt}
    \begin{subfigure}{0.49\textwidth}
        \includegraphics[width=\linewidth]{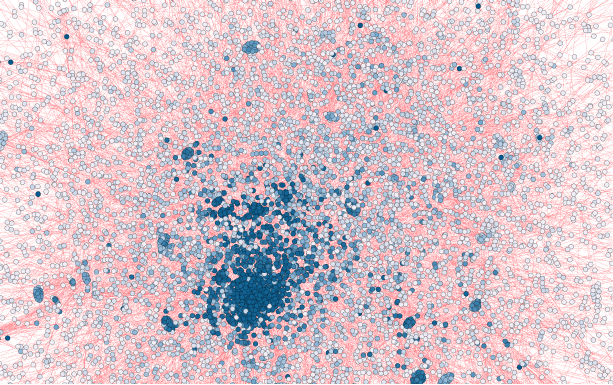}
        \caption{Four labels blocked}
        \label{fig:HepPh_final}
    \end{subfigure}
    \caption{The final spread on HepPh collaboration network with 12,008 nodes and 118,521 edges, (a) when no labels are blocked and (b) when four labels are blocked. 
    There are 50 seed nodes, 50 diffusion scenarios and 21 labels. The filling of the nodes depends on the number of scenarios each node is activated in. The darkest blue filling corresponds to the seed nodes, i.e., nodes which are activated in each of 50 scenarios. Arcs with blocked labels are removed from the graph in (b).} \label{fig:HepPh}
\end{figure}

An example of how blocking labels affects the final spread is displayed in Figure \ref{fig:HepPh} on a collaboration network instance whose details are explained in Section \ref{sec:results}.

\begin{proposition}
The objective function $f(K')$ of \eqref{eq:P} is neither submodular nor supermodular.
\end{proposition}

\begin{proof}
We show this by means of a counterexample given in Figure \ref{fig:example} where arc labels are shown next to the arcs. We assume that all the influence probabilities being equal to one, i.e., there is only one possible live-arc scenario, and all the arcs are live in that scenario. In this special case, it is enough to evaluate the number of nodes reachable from the seed set on the original graph to compute the expected spread.

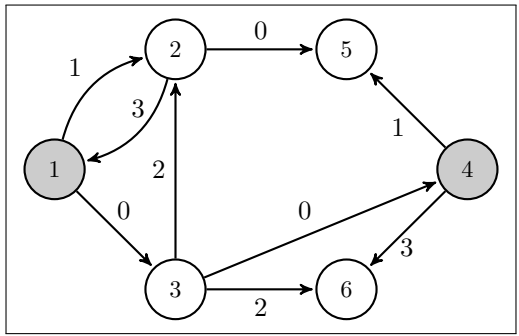
\begin{figure}[h!tb]
\centering
  \fbox{
    \begin{tikzpicture}[->,>=stealth',shorten >=1pt,auto,node distance=2.5cm,
                    thick, main node/.style={circle,draw,font=\sffamily\Large\bfseries}]
\tikzstyle{every state}=[fill=white,draw=black,text=black,scale=0.9]

  \node[state,fill=black!20!white] 		   (n1)             		  {1};
  \node[state]         (n2) [above right of=n1] 	  {2};
  \node[state]         (n3) [below right of=n1]       {3};
  \node[state]         (n5) [right of=n2] {5};
  \node[state]         (n0) [right of=n3] {6};
  \node[state,fill=black!20!white]         (n4) [below right of=n5] {4};
  \path (n1) edge[bend left]			  node {1} (n2)
            edge            node[xshift=-3pt] {0} (n3)    
        (n2) edge [bend left]			  node[above, yshift=2pt]  {3} (n1)
        	edge			  node {0} (n5)
        (n3)	edge	   node[left] {2} (n2)
         edge              node {0} (n4)
         edge              node[below] {2} (n0)
        (n4) edge 			  node[xshift=3pt] {1} (n5)
    	    edge			  node[below] {3} (n0)
        ;
\end{tikzpicture}
}
\caption{A small directed network with four labels.}\label{fig:example}
\end{figure}

Let $I=\{1,4\}$ be the seed set. For $f$ to be submodular, it needs to satisfy $f(X\cup \{i\})- f(X)\geq f(Y\cup \{i\})- f(Y)$ for all $X\subseteq Y \subseteq K$ and $i\in K\setminus Y$. For supermodularity $-f$ should be submodular. 
As a counterexample, consider the sets $X=\{3\}$ and $Y=\{2,3\}$. It is easy to compute $f(X)=6$ and $f(Y)=5$. If we choose $i=1$, then we get $f(X\cup \{i\})=f(\{3,1\})=6$ and $f(Y\cup \{i\})=f(\{2,3,1\})=3$, which violates the inequality we need for supermodularity. Now, suppose that we choose $i=6$ which yields $f(X\cup \{i\})=f(\{3,6\})=4$ and $f(Y\cup \{i\})=f(\{2,3,6\})=4$. These numbers violate the submodular inequality. We conclude that $f$ is neither submodular nor supermodular, as the requirements are not satisfied in the \MS{given counterexample}.
\end{proof}

\section{An arc-based compact formulation} \label{sec:arc-based}

In this section, we first present an arc-based compact formulation, and then show how \BD can be applied to this formulation. 

\subsection{The compact formulation}
%Recall that each arc in the directed graph has exactly one label, whereas it is possible to have \MS{parallel} arcs with different labels between a pair of nodes, and that blocking a label implies the removal of the arcs with the corresponding label from the network. The problem involves finding a set of labels to block within a given amount of budget so that the expected spread is minimized. 

For a scenario $\omega$, let $A^\omega_k$ denote the set of arcs with label $k$ in the live-arc graph. Let binary variable $x_k$, $k\in K$ be one iff label $k$ is blocked. Moreover, let binary variable $y^\omega_i$, $i \in V$, $\omega \in \Omega$ be one iff the spread of the harmful contagion reaches node $i$ in scenario $\omega$. Using this notation, a formulation of the \Problem can be given as follows, which we denote as \myref{ABF}.

\begin{align}
\mytag{ABF} \quad & \min \frac{1}{|\Omega|}\sum_{\omega\in \Omega} \sum_{i\in V} y_i^\omega \tag{ABF.OBJ}\label{mip00}\\
& \text{s.t.} \notag \\
& \hspace{15pt} \sum_{k \in K} c_k x_k \leq B \tag{ABF.BUD}\label{mip01}\\
& \hspace{15pt} y_i^\omega \geq 1 &i&\in I, \omega \in \Omega \tag{ABF.SEED} \label{mip03}\\
& \hspace{15pt}y_j^\omega\geq y_i^\omega -x_{k} &(&i,j) \in A_k^\omega , k\in K, \omega\in \Omega \tag{ABF.SPR} \label{mip02}\\
& \hspace{15pt}x_k\in \{0,1\} &k&\in K \notag\\
& \hspace{15pt}y_i^\omega\in \{0,1\} &i&\in V, \omega \in \Omega \notag
\end{align}

\RMS{
\begin{proposition}
The formulation \myref{ABF} models the \Problem.
\end{proposition}

\begin{proof}
The objective function \eqref{mip00} minimizes the sum of the nodes to which the harmful contagion spreads over all scenarios, which is the objective function of \Problem. Constraint \eqref{mip01} ensures that the budget for taking the measures is not exceeded. What now remains to show that the formulation correctly models \Problem is that for each scenario $\omega \in \Omega$ the spread of the harmful contagion is modeled correctly, considering the blocked labels indicated by variables $x_k$. Recall from definition \eqref{eq:spread} that the contagion spreads to a node $v \in V$ in a scenario $\omega$ (which means $y^\omega_v$ needs to have the value one), iff there exists a non-blocked path from a seed node $i \in I$ to $v$ in the live-arc graph with arc set $A^\omega$. Constraints \eqref{mip03} set the correct value of the $y$-variables for all the seed nodes (we model this with an inequalitiy instead of an equality to simplify the dualizing-step we need for the BD formulation below). Constraints \eqref{mip02} model the spread of the harmful contagion: if in a scenario $\omega$ the contagion has spread to node $i$ (i.e., $y^\omega_i$ is one) and there is an arc between $i$ and $j$ with label $k$ in the live-arc graph with arc set $A^\omega$, then the contagion will also spread to node $j$ (i.e., $y^\omega_j$ needs to have the value one to fulfill the constraint) unless label $k$ is blocked (i.e., $x_k$ is one). 
\end{proof}
}

We note that the number of variables in the formulation \myref{ABF} is in $O(|V||\Omega|)$ and the number of constraints is in $O(|A||\Omega|)$. Thus, while the formulation is compact, for large-scale graphs or a large number of scenarios, the size of the formulation can become prohibitive to be solved directly. To deal with this issue, we next present a \BD approach.

\subsection{Benders decomposition formulation}

It is easy to see that in formulation \myref{ABF} the $y$-variables can be relaxed to be continuous (i.e., $y^\omega_i\geq 0$), as long as the $x$-variables are binary. Due to the constraints \eqref{mip02} and the minimization-objective, the $y$-variables will also take binary values in a solution, \KT{for a fixed binary $x$}. This allows us to apply a scenario-based \BD to formulation \myref{ABF}. Let $\theta^\omega$, $\omega \in \Omega$ be a continuous variable to measure the spread of the harmful contagion in scenario $\omega$. The \BD reformulation of \myref{ABF} is given as follows.

\begin{align}
& \mytag{BEN} \quad \min\frac{1}{|\Omega|} \sum_{\omega\in \Omega} \theta^\omega \tag{BEN.OBJ} \label{eq:BMP_obj}\\
& \text{s.t.} \notag \\
& \hspace{15pt} \sum_{k \in K} c_k x_k \leq B \tag{BEN.BUD} \label{eq:BMP_c1}\\
& \hspace{15pt} \theta^\omega \geq \Theta^\omega(x) &\omega &\in \Omega \tag{BEN.OBJ2} \label{eq:BMP_c2}\\
& \hspace{15pt}x_k\in \{0,1\} &k&\in K \notag
\end{align}

Inequalities \eqref{eq:BMP_c2} provide a lower bound on the number of nodes which are reached by the spread of the harmful contagion for a given scenario $\omega \in \Omega$ and any $x \in \{0,1\}^K$. Benders optimality cuts \NA{can be derived} to model the function $\Theta^\omega(x)$ for each $\omega \in \Omega$ and a fixed $\bar x$, using linear programming (LP) duality. Note that for any fixed solution $\bar x$ of \myref{BEN}, the remaining subproblem in $y$ is always feasible, thus no Benders feasibility cuts are needed.  
Let $\alpha_i^\omega$ and $\beta_a^\omega$ be the dual variables associated with constraints \eqref{mip03} and \eqref{mip02}. Let $\delta^+_\omega(i)$ and $\delta^-_\omega(i)$, respectively, be the set of outgoing and incoming arcs of node $i$ in the live-arc graph of scenario $\omega$. We obtain the following dual for a particular scenario $\omega \in \Omega$.
\begin{align}
\mytag{DSABF} \quad & \max \sum_{i\in I} \alpha_i^\omega - \sum_{k\in K} \sum_{a\in A_k^\omega} \beta_a^\omega \bar{x}_{k} \label{eq:DSPr_obj} \tag{DSABF.OBJ}\\
& \text{s.t.} \notag \\
& \hspace{15pt} \alpha_i^\omega - \sum_{a\in \delta^{+}_\omega(i)} \beta_a^\omega   + \sum_{a\in \delta^{-}_\omega(i)} \beta_a^\omega \leq 1 &i& \in I \label{eq:DSPr_c1} \tag{DSABF.SEED} \\
& \hspace{15pt} - \sum_{a\in \delta^{+}_\omega(i)} \beta_a^\omega   + \sum_{a\in \delta^{-}_\omega(i)} \beta_a^\omega \leq 1 &i&\in V \setminus I \label{eq:DSPr_c2} \tag{DSABF.NONSEED}\\
& \hspace{15pt}\alpha_i^\omega\geq 0 &i&\in I \notag \\
& \hspace{15pt} \beta_a^\omega \geq 0 &a& \in A_k^\omega, k\in K \notag
\end{align}
%where $A_{i\omega}^{+}$ and $A_{i\omega}^{-}$ are, \NA{respectively,} the sets of outgoing and incoming arcs of node $i$ that are live in scenario $\omega$, respectively. 
For a solution $(\hat \alpha^\omega,\hat \beta^\omega)$ of \myref{DSABF}, the resulting Benders optimality cut is
%A single Benders optimality cut is
% \begin{equation}
% \theta \geq \sum_{r\in \Omega}\Big( \sum_{i\in I}  \hat{\alpha}_{ir} - \sum_{a\in A_r} \hat{\beta}_a^\omega x_{k(a)} 
% \Big)
% \end{equation} 
%or
\begin{equation}
\theta^\omega \geq \sum_{i\in I}  \hat{\alpha}_{i}^{\omega} - \sum_{k\in K}\sum_{a\in A_k^\omega} \hat{\beta}_a^\omega x_{k} \label{eq:BendersCut}\tag{BOC.ABF}
\end{equation} 

As \myref{DSABF} is a compact LP formulation, for the given scenario $\omega$ and the value of $\bar{x}$, \eqref{eq:BendersCut} can be separated by solving \myref{DSABF} as an LP. However, as an alternative method, the cuts can also be generated in a combinatorial and exact way, without solving the LP formulation of \myref{DSABF}. Let $G^\omega(\bar{x})$ denote the live-arc graph $G^\omega$ where the length of each arc $(i,j)\in A_k^\omega$ is equal to $\bar{x}_k$ for all $k\in K$. We also define $d_i^\omega(\bar{x})$ as the length of the shortest path from the seed set $I$ to $i$ on $G^\omega(\bar{x})$. In other words, $d_i^\omega(\bar{x})=\min_{j\in I} d_{ji}^\omega(\bar{x})$ where $d_{ji}^\omega(\bar{x})$ is the \RKT{length of the} shortest path from $j \in I$ to $i$ on $G^\omega(\bar{x})$.

\begin{definition}
    Given a solution $\bar{x}$ to \myref{ABF} and a scenario $\omega$, a node $i \in V$ is called \emph{reachable} \RKT{(from the seed set) }if $d_i^\omega(\bar{x})<1$ and not reachable otherwise. 
\end{definition}

\begin{definition}
    The path that is associated with the shortest path distance $d_i^\omega(\bar{x})$ from the seed set to the reachable node $i$ is the \emph{activation path} of $i$ and each reachable node has exactly one activation path (ties are broken arbitrarily). 
    \end{definition}

\begin{example}

\RKT{Recall the instance in Figure \ref{fig:problem} where the seed set is $I=\{1,4\}$ and the blocking decision is \RMS{$\bar{x}=(0,0,1,0)$}. The length of an arc with label $k$ is calculated as $\bar{x}_k$ and the resulting shortest path distances to each node from $I$ in the first scenario are $d_1(\bar{x})=d_3(\bar{x})=d_4(\bar{x})=d_5(\bar{x})=d_6(\bar{x})=0$ and $d_2(\bar{x})=1$ (see Figure
 \ref{fig:activation_path}). Thus, all nodes are \emph{reachable} except node 2.
 %which is the only node that is not activated.
 The \emph{activation paths} of 1 and 4 start and end at themselves as they are seed nodes (i.e., they are active initially), whereas the activation paths of nodes 3, 5, and 6 are 1-3,  4-5, and 4-6 respectively. Notice that there are two paths to each of nodes 5 and 6, one from node 1 and the other from node 4, and only the paths from node 4 have a length of less than one. The reason is that there are no blocked arcs on those paths, so node 4 can reach and activate them.

}

 \begin{figure}[h!tb]
\centering
 \fbox{
    \begin{tikzpicture}[->,>=stealth',shorten >=1pt,auto,node distance=2.5cm,
                    thick, main node/.style={circle,draw,font=\sffamily\Large\bfseries}]
\tikzstyle{every state}=[fill=white,draw=black,text=black,scale=0.9]

  \node[state,fill=black!20!white,label=above :{\footnotesize $d_1(\bar{x})=0$}]  (n1)   {1};
  \node[state,label=above :{\footnotesize $d_2(\bar{x})=1$} ]         (n2) [above right of=n1] 	  {2};
  \node[state,fill=black!50!white,label=below :{\footnotesize $d_3(\bar{x})=0$}]         (n3) [below right of=n1]       {3};
  \node[state,fill=black!50!white,label=above :{\footnotesize $d_5(\bar{x})=0$}]         (n5) [right of=n2] {5};
  \node[state,fill=black!50!white,label=below :{\footnotesize $d_6(\bar{x})=0$}]         (n0) [right of=n3] {6};
  \node[state,fill=black!20!white,label= above:{\footnotesize $d_4(\bar{x})=0$}]         (n4) [below right of=n5] {4};
  \path (n1)  %edge[bend left]			  node {1} (n2)
        edge            node[xshift=-3pt] {0} (n3)    
        (n2) %edge [bend left]			  node[above, yshift=2pt]  {3} (n1)
        	edge  			  node {0} (n5)
        (n3)	edge [dashed]	   node[left] {2} (n2)
         %edge              node {0} (n4)
         edge     [dashed]         node[below] {2} (n0)
        (n4) edge 			  node[xshift=3pt] {1} (n5)
    	    edge			  node[below] {3} (n0)
        ;
\end{tikzpicture}
}
\caption{Reachable nodes and activation paths of a scenario for a given blocking decision }
\label{fig:activation_path}
 \end{figure}
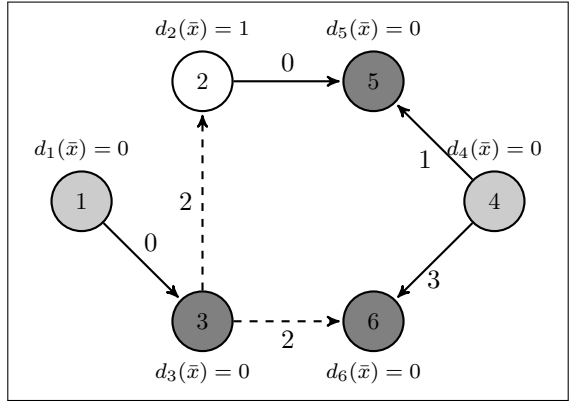

  \end{example}

\begin{theorem}
Consider the solution \RKT{$(\hat \alpha^\omega,\hat \beta^\omega)$} to \myref{DSABF} where \RKT{$\hat{\alpha}_{i}^{\omega}$ is
%the number of nodes that are reachable via a path emerging from node $i\in I$, and $\hat{\beta}_a^\omega$} is the number of nodes whose activation path contains the arc $a \in A^\omega$
the number of activation paths starting at node $i\in I$ and $\hat{\beta}_a^\omega$ is the number of activation paths passing through arc $a \in A^\omega$.} This solution \RKT{$(\hat \alpha^\omega,\hat \beta^\omega)$} is an optimal solution to \myref{DSABF}. 
\label{prop:DSP_solution}
\end{theorem}

\begin{proof}
See \RMS{Appendix} \ref{appendix_proof}.
\end{proof}

We note that the optimal solution \RKT{$(\hat \alpha^\omega,\hat \beta^\omega)$} satisfying the condition in Theorem \ref{prop:DSP_solution} can be computed using a shortest path algorithm. In the cut obtained, the constant part $\sum_{i\in I} \hat{\alpha}_{i}^{\omega}$ corresponds to the total number of nodes that are fully ($y_i^\omega=1$) or partially ($0<y_i^\omega<1$) activated by $I$ in scenario $\omega$, \NA{when} $x=\bar{x}$. The coefficients $\sum_{a\in A_k^\omega} \hat{\beta}_a^\omega$ of each $x_k$ correspond to the maximum possible reduction in the number of 
%newly activated
active nodes if label $k$ is blocked. Being able to solve the subproblems to obtain Benders optimality cuts in a combinatorial fashion by shortest path computations instead of solving LPs can be very useful. However, as we show in the next section, it is possible to obtain a different formulation, which also allows for such combinatorial computation of Benders optimality cuts, and additionally allows for lifting of these cuts. The drawback of the formulation of the next section is that it has exponentially many constraints. Thus, compared to the compact formulation provided in this section, it is not possible to give it directly to an \ILP solver\MS{\ if one wants to solve the problem without implementing a BD,} \KT{or to solve the Benders subproblems as LPs. We describe how to handle this drawback in the following section.}

\section{A path-based formulation} \label{sec:path-based}

In this section, we first describe the path-based formulation, and then show that \RKT{the LP-relaxations of the path-based formulation and the arc-based formulation introduced} in the previous section have the same optimal objective value. Next, we present a pair of valid inequalities which are obtained by down-lifting the constraints of the path-based formulation and show that these inequalities \MS{can} improve the LP-relaxation. We then discuss how to apply \BD to the path-based formulation and show that despite the exponential number of inequalities in the formulation, the Benders optimality cuts (based on the non-lifted inequalities) can be obtained in a combinatorial fashion in polynomial time. Finally, we discuss why the separation of Benders optimality cuts is not straightforward in the case of the path-based formulation with the down-lifted version of the constraints. 

\subsection{Formulation}

The path-based formulation is based on directly modeling the fact that the contagion reaches a node $i \in V$ in a scenario $\omega$ iff there exist an unblocked path from the seed set $I$ %$S$
to this node $i$ in the live-arc graph of the scenario (see equation \eqref{eq:spread} in Section \ref{sec:diff} where this is stated in a cumulative fashion). Let $P^\omega(i)$ denote the set of all paths from the seed set to a given node $i$ in scenario $\omega$. For the ease of notation, we assume that for a seed node $i \in I$ this set contains an empty path (which by definition is unblockable since there is no arc on it).
%We can define an alternative model by summing up constraints \eqref{mip02} for all possible paths from the seed set $I$ to any node $i \in V^\prime$. Let $P^\omega(i)$ denote the set of such paths for a given node $i$ in scenario $\omega$. 
%
Moreover, for a given path $p \in P^\omega(i)$, let $K_p$ be the set of labels that is encountered on $p$ and $h_p(k)$ be the number of times label $k$ occurs on this path, i.e., the number of arcs $(i,j) \in p$ with label $k$. 
%\footnote{\KT{Do we need a different notation than $d$ as we use it for shortest path distances in Section 3.2} \NA{I think so, we can use $h_p(k)$}}
%
We obtain the following formulation:
%For ease of notation let $P^\omega(i)=\emptyset$ for $i \in I$.

\begin{align}
\mytag{PBF} \quad& \min \frac{1}{|\Omega|}\sum_{\omega\in \Omega} \sum_{i\in V} y_i^\omega \tag{PBF.OBJ} \label{nmip00}\\
& \text{s.t.} \notag \\
& \hspace{15pt} \sum_{k \in K} c_k x_k \leq B \label{nmip01} \tag{PBF.BUD}\\
& \hspace{15pt}y_i^\omega\geq 1 -\sum_{k \in K_p}h_p(k) x_{k} &i&  \in V,  \omega\in \Omega , p \in P^\omega(i) \label{nmip02} \tag{PBF.PATHS} \\
& \hspace{15pt}x_k\in \{0,1\} &k&\in K \notag \\
& \hspace{15pt}y_i^\omega \in \{0,1\} &i&\in V, \omega \in \Omega \notag
\end{align}

\RMS{
\begin{proposition}
    The formulation \myref{ABF} models the \Problem.
\end{proposition}
\begin{proof}
The objective function \eqref{nmip00} and the budget constraint \eqref{nmip01} are the same as in the formulation \myref{ABF}. Constraints \eqref{nmip02} directly model the spread for each scenario $\omega$ following its definition \eqref{eq:spread}: they ensure that for a given $\bar x$ the variable $y_i^\omega$ needs to have the value one if there exists at least one unblocked path from the seed set.
\end{proof}
}
We note that there can be an exponential number of these inequalities \eqref{nmip02} since there can be an exponential number of paths.

Next, we compare the LP-relaxations of \myref{ABF} and \myref{PBF}. In these relaxations, the binary variables $x_k$, $k\in K$ and $y_i^\omega$, $i\in V, \omega \in \Omega$ are relaxed to $0\leq x_k\leq 1$ and $0\leq y_i^\omega\leq 1$.

\begin{theorem}
The optimal objective function values of the LP-relaxations of \myref{ABF} and \myref{PBF} are the same.
\end{theorem}

\begin{proof}
We note that constraints \eqref{nmip02} of the formulation \myref{PBF} can be obtained by summing up constraints \eqref{mip03} and \eqref{mip02} for all possible paths. Aside from these constraints, both formulations are the same. Thus, the feasible region of the LP-relaxation of the formulation \myref{ABF} is at most as large as the feasible region of the LP-relaxation of formulation \myref{PBF}. Consequently, the optimal objective value of the LP-relaxation of the formulation \myref{ABF} is always greater than or equal to the optimal objective value of the LP-relaxation of the formulation \myref{PBF}.

To show equality, we next show that the optimal solution $(\hat x, \hat y)$ of the LP-relaxation of the formulation \myref{PBF} is also feasible for the LP-relaxation of the formulation \myref{ABF}. Suppose this is not the case. This means $(\hat x, \hat y)$ either violates one of the constraints \eqref{mip03} or one of the constraints \eqref{mip02}. We first note that $(\hat x, \hat y)$ cannot violate any of the constraints \eqref{mip03} as they actually are a special case of constraints \eqref{nmip02} for $i \in I$. Thus, suppose $(\hat x, \hat y)$ violates a constraint \eqref{mip02} which would mean $\hat y^\omega_j<\hat y^\omega_i-\hat x_k$ for some $\omega \in \Omega$, $k \in K$, $(i,j) \in A^\omega_k$. Due to $\hat y^\omega_j \geq 0$ this is only possible if $y^\omega_i>0$.  Due to the minimization objective and constraints \eqref{nmip02}, there must be a path $p^* \in P^\omega(i)$ with $\hat y^\omega_i =1 -\sum_{k' \in K_{p^*}}h_{p^*}(k') \hat x_{k'}$. Moreover, the path $q^* \in P^\omega(j)=p^* \cup (i,j)$ must also exist as otherwise there would be no constraint \eqref{mip02} for the considered $\omega \in \Omega$, $k \in K,(i,j) \in A^\omega_k$. For this path it must hold that $\hat y^\omega_j \geq 1 -\sum_{k' \in K_{q^*}}h_{q^*}(k') \hat x_{k'}$ due to constraints \eqref{nmip02}. From this we get $\hat y^\omega_j\geq \hat y^\omega_i-\hat x_k$, which is a contradiction to our assumption $\hat y^\omega_j<\hat y^\omega_i-\hat x_k$. As a consequence the solution $(\hat x, \hat y)$ is feasible, \NA{and thus optimal,} for the LP-relaxation of the formulation \myref{ABF}.

%\footnote{\KT{if the arc is $(i,j)$, $q^*$ should be the longer path}\NA{you are right, what about the switch of $i$ and $j$?}}

\end{proof}

\RKT{Although the LP-relaxations of both \myref{ABF} and \myref{PBF} provide the same lower bounds, constraints \eqref{nmip02} can be down-lifted as described in the following proposition, which could decrease the lower bound obtained from the LP-relaxation of \myref{PBF}.
%Although the strengths of the LP-relaxations of \myref{ABF} and \myref{PBF} are the same, constraints \eqref{nmip02} can be down-lifted as described in the following \RKT{proposition}, which could improve the LP-relaxation of \myref{PBF}.
}

\begin{proposition}
Given a scenario $\omega \in \Omega$, node $i \in V$ and path $p \in P^\omega(i)$, inequalities
\begin{equation}
    y_i^\omega\geq 1 -\sum_{k \in K_p} x_{k} \tag{PBF.PATHSL}
    \label{nmip02lifted}
\end{equation}
are valid.
\end{proposition}

\begin{proof}
   For each $\omega \in \Omega$, node $i \in V$ and path $p \in P^\omega(i)$ a feasible solution $(\hat x, \hat y)$ must \NA{satisfy the associated} constraint~\eqref{nmip02}. The right-hand-side of this constraint is one iff all variables $\hat x_k$ for $k \in K_p$ are zero, otherwise the right-hand-side of this constraint is zero or smaller than zero. As each $y^\omega_i$ can only take values zero or one, this means the constraint~\eqref{nmip02} becomes redundant as soon as one of the variables in $\hat x_k$ for $k \in K_p$ has the value one. It is easy to see that constraint~\eqref{nmip02lifted} behaves exactly the same.
\end{proof}

Let $\mytag{PBF.L}$ denote the lifted path-based formulation where inequalities \eqref{nmip02lifted} replace inequalities \eqref{nmip02}.
\begin{proposition}
There exist instances, where the optimal objective value of the LP-relaxation of \myref{PBF.L} is strictly greater than the optimal objective value of the LP-relaxation of \myref{PBF}.
\end{proposition}

\begin{proof}
Consider an instance with $V=\{a,b,c\}$, $I=\{a\}$, $|\Omega|=1$ and the following arcs exist in the live-arc graph for the single scenario: $(a,b), (b,c)$. Moreover, we have a single label $\ell$ (both arcs have the same label), and blocking this label has a cost of two units. The budget $B$ is one. Thus, we cannot block the label $\ell$, and in the optimal solution, the contagion spreads to all three nodes. Hence the optimal objective value is equal to three. In the LP-relaxation of \myref{PBF}, we get an optimal value of $1.5$, since $\hat x_\ell=0.5$ is feasible, and we get $\hat y^\omega_a=1$, $\hat y^\omega_b=0.5$ and $\hat y^\omega_c=0$. In the LP-relaxation of \myref{PBF.L}, we get an optimal value of two since for $\hat x_\ell=0.5$ we get $\hat y^\omega_a=1$, $\hat y^\omega_b=0.5$ and $\hat y^\omega_c=0.5$. The difference in the value of $\hat y^\omega_c$ in both formulations is caused by the fact that in \myref{PBF}, we have $-2x_k$ on the right-hand-side of the inequality \eqref{nmip02} for the path to $c$, while in \myref{PBF.L}, we have $-x_k$ on the right-hand-side of the inequality \eqref{nmip02lifted}.
\end{proof}

\subsection{Benders decomposition formulations}
\label{sec:nonlifted_benders}

We first focus on \BD for the formulation \myref{PBF} and then extend it to \BD for the formulation \myref{PBF.L}. As a starting point, we observe that similar to the formulation \myref{ABF}, we can relax the $y$-variables to $y^\omega_i\geq0$ in \myref{PBF} (and also \myref{PBF.L}), and still get binary $y$-values in an optimal solution for a fixed $\bar x$. The Benders master problem is the same as \myref{BEN} and only the subproblem for a fixed $\bar x$ changes with the path-based formulations. For a given $\bar{x}$, the number of nodes \NA{to which} the contagion has spread in scenario $\omega$ can be obtained via solving the following (primal) subproblem when considering \myref{PBF}.
\begin{align}
\mytag{PSPBF}\quad& \min \sum_{i\in V} y_i^\omega \label{eq:nSPr_obj} \tag{PSPBF.OBJ}\\
& \text{s.t.} \notag \\
& \hspace{15pt}y_i^\omega\geq 1 -\sum_{k \in K_p}h_p(k) \bar{x}_{k} & & i \in V, p \in P^\omega(i) \label{eq:nSPr_c2} \tag{PSPBF.PATHS}\\
%& \hspace{15pt}y_i^\omega\geq 1 &i&\in I\\
& \hspace{15pt}y_i^\omega\geq 0 &i&\in V \notag
\end{align}

Let $\gamma^\omega_{ip}$ be dual variables associated with constraints \eqref{eq:nSPr_c2}. The dual of the \myref{PSPBF} is given as follows.

\begin{align}
\mytag{DSPBF}\quad& \max \sum_{i\in V} \sum_{p\in P^\omega(i)} \big ( 1-\sum_{k \in K_p}h_p(k) \bar x_{k} \big ) \gamma^\omega_{ip} \label{eq:nDSPr_obj} \tag{DSPBF.OBJ}\\
& \text{s.t.} \notag \\
& \sum_{p\in P^\omega(i)} \gamma^\omega_{ip} \leq 1 & & i \in V \tag{DSPBF.PATHS} \label{eq:nDSPr_c2}\\
& \hspace{15pt}\gamma_{ip} \geq 0 &i&\in V, p\in P^\omega(i) \notag
\end{align}

For a solution $\hat \gamma^\omega$, the resulting Benders optimality cut is

\begin{equation}
   \theta^\omega \geq \sum_{i\in V} \sum_{p\in P^\omega(i)} \hat{\gamma}^{\omega}_{ip}   -\sum_{i\in V} \sum_{p\in P^\omega(i)} \sum_{k \in K_p} h_p(k) \hat{\gamma}^{\omega}_{ip} x_{k} \label{eq:BendersCut2}. \tag{BOC.PBF}
\end{equation}

The formulation \myref{DSPBF} includes an exponential number of variables and solving it as an LP requires the enumeration of all $p\in P^\omega(i)$ for all $i\in V \setminus I$ (for $i \in I$, the empty path is trivially the best and thus the only one we need to include). Therefore, \myref{PBF} is not a suitable formulation for a \BD method where the subproblems are solved as LPs. However, the structure of \myref{DSPBF} allows obtaining the cut \eqref{eq:BendersCut2} in polynomial time in a combinatorial way, using the characterization of an optimal solution to \myref{DSPBF} which is described in the following theorem. 

\begin{theorem}
Let $\hat{\gamma}^\omega$ be a solution to \myref{DSPBF} such that $\hat \gamma^\omega_{ip}$ is equal to one if $i$ is reachable and path $p$ is the activation path of $i$, and equal to zero otherwise.  Then, $\hat{\gamma}^\omega$ is an optimal solution to \myref{DSPBF}. 
\label{prop:BendersCut_Path}
\end{theorem}

\begin{proof}
As can be readily seen, \myref{DSPBF} decomposes into $|V|$ subproblems for the nodes $i \in V$, where each has the objective $\max \sum_{p\in P^\omega(i)} \big ( 1-\sum_{k \in K_p}h_p(k) \bar x_{k} \big ) \gamma^\omega_{ip}$. The value of $\sum_{k \in K_p}h_p(k) \bar x_{k}$ is equal to the length of the path $p$ on $G^\omega(\bar{x})$. For given $i$, let $p^\prime \in P^\omega(i)$ be the path with the shortest length on $G^\omega(\bar{x})$ (ties broken arbitrarily). If $p^\prime$ has a length less than one, i.e., $i$ is reachable and $p^\prime$ is its activation path, setting $\gamma^\omega_{ip^\prime}=1$ and $\gamma^\omega_{ip}=0$ for all $p\neq p^\prime$ gives an optimal decision for the problem associated with $i$ as this results in the largest possible contribution to the objective function.  If $p^\prime$ has a length of at least one, then the objective function coefficients of all $\gamma^\omega_{ip}$ are non-positive and thus we need to have $\gamma^\omega_{ip}=0$ for all $p\in P^\omega(i)$ in an optimal solution. 
\end{proof}

For an optimal value of $\hat{\gamma}^\omega$ as described in Theorem \ref{prop:BendersCut_Path}, $\sum_{i\in V} \sum_{p\in P^\omega(i)} \hat{\gamma}^\omega_{ip}$ gives the number of nodes that are reachable in scenario $\omega$. The coefficient of $x_k$ in \eqref{eq:BendersCut2} represents the number of times an arc with label $k$ appears on the activation paths of all nodes, as in \eqref{eq:BendersCut}.

\begin{remark}
It is easy to observe that $d_i^\omega(\bar{x})\in \mathbb{N}$ for $\bar{x}\in \{0,1\}^{|K|}$. In this case, if the shortest path to node $i$ includes at least one blocked arc, i.e., some arc $(i,j)\in A_k^\omega$ such that $\bar{x}_k=1$, then $d_i^\omega(\bar{x})\geq 1$ and node $i$ is not reachable. Therefore, instead of the shortest path lengths, it is enough to have the information 
if node $i$ is accessible from the seed set $I$ via the arcs in $\cup_{k|\bar{x}_k=0}A_k^\omega$. By Theorem \ref{prop:BendersCut_Path}, $\hat{\gamma}^\omega_{ip}=0$ 
%for any $a\in \cup_{k|\bar{x}_k=1}A_k^\omega$. 
for all paths $p$ such that $\sum_{k\in K_p} \bar{x}_k\geq 1$.
As a result, we get 
that the optimal solution value of \myref{DSPBF} in this case is
$\sum_{i\in V} \sum_{p\in P^\omega(i)} \hat{\gamma}^\omega_{ip}$. This is equal to the number of nodes \NA{to which} the contagion can spread in scenario $\omega$ for given $\bar{x}$, since a reachable node for a given $\bar{x}$ corresponds to an activated node ($y_i=1$) in the diffusion process.
\label{remark:integerSep}
\end{remark}

%\KT{TODO: rewrite the following separation paragraphs}

Based on Theorem \ref{prop:BendersCut_Path} and Remark \ref{remark:integerSep}, we can make use of a simple graph search algorithm to separate the cut \eqref{eq:BendersCut2} via computing an optimal dual subproblem solution $\hat{\gamma}^\omega$, when $\bar{x}\in\{0,1\}^{|K|}$. 
In this method, for a given scenario $\omega\in \Omega$, first the nodes that are accessible from $I$ via paths consisting only of the arcs in $\cup_{k:\bar{x}_k=0} A_k^\omega$ are determined. 
There is a violated cut if the number of accessible nodes is strictly greater than $\bar{\theta}^\omega$. In this case, we compute the cut coefficients $\sum_{i\in V} \sum_{p\in P^\omega(i)} h_p(k) \hat \gamma^\omega_{ip}$ for each $k \in K$ by tracking along the activation paths. However, using a single source shortest path algorithm, an exact polynomial time separation of the cuts is also possible for any fractional solution $\bar{x}$. We discuss this separation algorithm in detail in Section \ref{sec:separation}.

For the \BD of \myref{PBF.L}, we observe that the resulting dual subproblem, denoted as $\mytag{DSPBF.L}$, is almost the same as \myref{DSPBF} in which only the objective function changes to

\begin{equation}
    \max \sum_{i\in V} \sum_{p\in P^\omega(i)} (1-\sum_{k \in K_p} \bar x_{k}) \gamma^\omega_{ip} \tag{DSPBF.L.OBJ} \label{eq:nDSPr_obj_lifted}
\end{equation}

and consequently, the Benders optimality cut changes to
\begin{equation}
   \theta^\omega \geq \sum_{i\in V} \sum_{p\in P^\omega(i)} \hat{\gamma}^\omega_{ip}   -\sum_{i\in V} \sum_{p\in P^\omega(i)} \sum_{k \in K_p} \hat{\gamma}^\omega_{ip} x_{k}. \label{eq:BendersCut3} \tag{BOC.PBF.L}
\end{equation}

Note that, differently from \eqref{eq:BendersCut} and \eqref{eq:BendersCut2}, the coefficient of $x_k$ in \eqref{eq:BendersCut3} is equal to the number of activation paths that include at least one arc with label $k$ (and not the cumulative number of occurences an arc with label $k$ appears on the activation paths).

Similar to \myref{DSPBF}, \myref{DSPBF.L} has an objective function that is decomposable in $i$. However, an optimal solution to \myref{DSPBF.L} cannot be found via computing shortest paths on $G^\omega(\bar{x})$.
%as was the case with \DSPp. 
Instead, it is required to find the path $p\in P^\omega(i)$ with the smallest value of $\sum_{k\in K_p}\bar{x}_k$, for each $i\in V$. 
Note that for binary valued $\bar{x}$, $\sum_{k\in K_p}\bar{x}_k=0$ (for an optimal path) if $d_i^\omega(\bar{x})=0$, and $\sum_{k\in K_p}\bar{x}_k\geq 1$ otherwise. Therefore, for binary valued $\bar{x}$, the graph search approach still works for finding a path minimizing $\sum_{k\in K_p}\bar{x}_k$. This means for binary valued $\bar{x}$, we can solve the separation problem for \eqref{eq:BendersCut3} exactly in polynomial time.
However, for fractional $ \bar{x}$, a shortest path approach is likely to fail because once we determine the shortest path to a node $i$, it does not necessarily lead to the best path(s) via $i$ for its successors. For example, suppose node $i$ can be reached by two different labels, say $0$ and $1$ from the seed set $I$, with $\bar x_0=0.4$ and $\bar x_1=0.5$. Thus, the shortest path to $i$ is achieved by taking the arc with label $0$. However, suppose now there is also a node $i'$, which is reachable from $i$ with an arc, which has label $1$ again. The path to $i'$ continuing from the shortest path to $i$ has the value $\bar x_0+\bar x_1=0.9$, while the path to $i'$ taking only label $1$ has the value $0.5$ and thus is better. Therefore, for solving the separation problem for \eqref{eq:BendersCut3} in case $\bar{x}$ is fractional we propose a heuristic method which is explained in Section \ref{sec:separation}. 

\section{Algorithmic details} \label{sec:algorithmic_details}

We propose solution frameworks for solving the Benders master problem \myref{BEN} within a \BC algorithm and add the Benders optimality cuts on-the-fly as they are needed. We develop such a \BD algorithm based on the arc-based formulation \myref{ABF}, and another one based on the path-based formulation \myref{PBF} and its lifted version \myref{PBF.L}. In the \BD algorithm based on the arc-based formulation, which serves as a baseline in our computational experiments, the subproblem for obtaining the Benders optimality cuts is compact and can be solved via linear programming. %Thus, the resulting \BD is a standard ``textbook'' implementation of the \BD algorithm.
The one based on the path-based formulation is the main algorithm considered in this work where we generate the Benders optimality cuts \eqref{eq:BendersCut2} (resp., \eqref{eq:BendersCut3}) with combinatorial algorithms, which we describe in Section \ref{sec:separation}. Moreover, this \BD algorithm also contains additional speed-ups such as scenario-dependent initial seed sets, cut sampling, a starting heuristic and initial cuts. All these enhancements are also discussed in this section. In our computational study, we do not only compare the performance of this enhanced \BD algorithm against the standard implementation of the \BD algorithm based on \myref{ABF}, but also analyze the effect of the individual enhancements.

\subsection{Separation algorithms for Benders optimality cuts}
\label{sec:separation}

Let $(\overline{\mathrm{BEN}})$
%\footnote{\KT{the overline looks too long to me, but I couldn't shorten it while still using a ref \NA{why don't we write simply $(\overline{\mathrm{BEN}})$ without using the reference?}}}
denote the problem to be solved at a node of the branch-and-cut tree. The problem $(\overline{\mathrm{BEN}})$ is a restricted version of the LP-relaxation of \myref{BEN}, excluding the constraints \eqref{eq:BMP_c2}, but (potentially) containing previously obtained Benders optimality cuts and branching decisions. Let $(\bar{\theta},\bar{x})$ be a feasible solution to $\overline{\myref{BEN}}$. If the optimal objective value of the considered dual subproblem (i.e., \myref{DSABF}, \myref{DSPBF} or \myref{DSPBF.L} depending on the chosen algorithmic setting) is greater than $\bar{\theta}^\omega$ for any $\omega\in \Omega$, then a violated Benders optimality cut exists and $(\bar{\theta},\bar{x})$ needs to be cut off. For the correctness of the algorithm, it is only required to check the existence of a violated cut for the solutions that satisfy the binary constraint $x\in\{0,1\}^{|K|}$, as fractional solutions are handled via branching in any case. However, separating fractional solutions can tighten the dual bounds and therefore improve the efficiency of the \BC algorithm. Thus, we also consider separating fractional solutions in some of our algorithmic settings which are considered in our computational study in Section \ref{sec:results}.

%If a solution $(\bar{\theta},\bar{x})$ satisfies the binary constraints but is not feasible for constraints \eqref{eq:BMP_c2}, at least one violated Benders optimality cut must be added to \myref{BEN}.

\paragraph{A basic separation algorithm}

As discussed above, in the basic implementation of the \BD algorithm, we solve \myref{DSABF}, which is a compact formulation, for each scenario $\omega \in \Omega$ as an LP and use its optimal solution to obtain the Benders optimality cut \eqref{eq:BendersCut}.

\paragraph{A combinatorial separation algorithm}

In order to solve the subproblem \myref{DSPBF} in a combinatorial way for given $\omega$ and $\bar{x}$ to obtain Benders optimality cuts \eqref{eq:BendersCut2}, we use the following algorithm. It works for both binary and fractional $\bar{x}$. As explained in the proof of Theorem \ref{prop:BendersCut_Path}, we need to compute the shortest path distances $d_i^\omega(\bar{x})$ for each $i\in V$. Note that for all seed nodes $i \in I$, this distance is trivially zero. Thus, we only need to focus on the nodes $V'=V \setminus I$. To do so, we use the priority queue version of the Dijkstra's shortest path algorithm, where we initialize the node distances as one for $i\in V'$ and as zero for $i\in I$, since $d_i^\omega(\bar{x})$ has to be strictly less than one for node $i$ to be \emph{reachable}. Furthermore, we stop the algorithm when the minimum distance in the course of the algorithm reaches one, as it is not possible to find any reachable nodes (i.e., nodes with $d_i^\omega(\bar{x})<1$) after that point. To compute the value of the cut coefficients (i.e., the value of the dual variables $\gamma^\omega_{ip}$), we backtrack on the paths after the Disjktra's algorithm terminates. 
%Recall that we need the  $\hat{\beta}$ values to decide if there is a violated cut \eqref{eq:BendersCut2} at a fractional solution $\bar{x}$. 
The overall separation procedure is described in more detail in Algorithm \ref{alg:Separation_scenario} which also has an option to heuristically solve the  subproblem \myref{DSPBF.L}, i.e., to heuristically obtain the lifted Benders optimality cuts \eqref{eq:BendersCut3}. This heuristic modification is discussed later below. 

In Algorithm \ref{alg:Separation_scenario} the option $\texttt{lift}=\texttt{N}$ means that the heuristic-part is not invoked. For a given scenario $\omega$, Algorithm \ref{alg:Separation_scenario} starts with initializing a priority queue \texttt{PQ} with the seed nodes. For each node $i\in V^\prime$, we denote by $d_i$ the distance from the seed set,  and by $K_{p_i}$ the set of labels used on the path $p_i$ to node $i$. 
We also initialize the variable $d_{p_i}(k)$, the number of times label $k$ is encountered on $p_i$, with zero.
If the new distance is smaller than the current one, it is updated as well as the path $K_{p_i}$ and $d_{p_i}(k)$. Once the priority queue is empty or the minimum distance is already one, we proceed to the last part where we compute the cut coefficients \KT{$c^\omega_k$, $\forall k,$ and the constant part of the cut $C^\omega$}. In line 12, we add \MS{a small} disturbance factor $\epsilon>0$ to each newly calculated distance so that we avoid using unnecessarily long paths \MS{(i.e., containing many arcs with labels $k$ where $\bar x_k=0$)} which could lead to weaker cuts due to double counting.

\paragraph{A modified combinatorial separation to \MS{heuristically} obtain lifted optimality cuts}

In order to obtain a lifted cut \eqref{eq:BendersCut3}\sloppy, we propose two heuristic \MS{enhancements which are} obtained by modifying the original separation algorithm in the distance calculation and path determination steps. In the first \MS{enhancement} which we call the \emph{posterior lifting}, we calculate the node distances using $\bar{x}$ as before, and while determining the number of times that a label is used on a path $p_i$ we count each label at most once. In other words, the maximum value that $d_{p_i}(k)$ can take is one (in line 18 of Algorithm \ref{alg:Separation_scenario}). This lifting option will be indicated by $\texttt{lift=P}$. 
In the second version which we refer to as \emph{heuristic lifting}, we change the way we compute the path lengths in a way that can help finding a path $p\in P^\omega(i)$ with a smaller value of $\sum_{k\in K_p}\bar{x}_k$, for each $i\in V^\prime$. To this end, when updating the distance of a node $v$ through another node $u$, we treat as if the length of arc $(u,v)$ is zero if the label of $(u,v)$ is already visited on the path to $u$, independent of the value of $\bar{x}_k$ where $(u,v)\in A^\omega_k$. Thus, we provide advantage to the paths that use a smaller number of distinct labels. As in the case of posterior lifting, we count each label at most once on a path, so that we obtain the lifted cut \eqref{eq:BendersCut3}, not \eqref{eq:BendersCut2}. We denote this lifting approach by setting the parameter $\texttt{lift=H}$. Note that, as can be seen in line 14, under this setting we do not punish the usage of long paths with the $\epsilon$ value since the length of the path does not have a direct effect on the strength of \eqref{eq:BendersCut3}.

\begin{algorithm}[htbp]
\caption{\texttt{Separation($\omega, \bar{\theta}, \bar{x}$,\texttt{lift})}}\label{alg:Separation_scenario}
\SetKwInOut{Input}{Input}\SetKwInOut{Output}{Output}
\Input{An infeasible master solution ($\bar{\theta},\bar{x}$)}
\Output{A (possibly violated) cut $\theta^\omega\geq C^\omega - \sum_{k\in K}c^\omega_k x_k$} 
Initialize a priority queue $\texttt{PQ}\leftarrow \emptyset$\; 
Initialize $C^\omega\leftarrow 0$ and $c_k^\omega\leftarrow 0 $, $\forall k \in K$\;
%Initialize $\texttt{reached}\leftarrow \texttt{false}$\;
Initialize the distance $d_i\leftarrow 1$, set of path labels $K_{p_i} \leftarrow \emptyset$, and label counts $d_{p_i}(k)\leftarrow 0$, $\forall i\in V, k\in K$\;
\For{each $i\in I$}
{
Set $d_i=0$, insert $(d_i,i)$ into \texttt{PQ}\;
}
\While{\texttt{PQ} $\neq \emptyset$}
{
    Choose a minimum distance unreached node $u \in \texttt{PQ}$ , label $u$ as \texttt{reached}\;
    \If{$d_u\geq 1$}{\textbf{break}\;}
    
\For{each outgoing arc of $(u,v)\in A^\omega$}
{
\If{$v$ is not \texttt{reached}}
{
Determine the label $k$ of arc $(u,v)\in A_k^\omega$, set $\hat{d}_v \leftarrow d_u + \max\{\bar{x}_k, \epsilon \}$ \;
\If{\texttt{lift=H} and $k\notin K_{p_v}$}
{
Set $\hat{d}_v \leftarrow d_u$ \;
}
\If{$\hat{d}_v < d_v$}
{
Set $d_v \leftarrow \hat{d}_v$, $K_{p_v} \leftarrow K_{p_u} \cup \{k\}$ \;
\If{\texttt{lift}=\texttt{N} or $d_{p_v}(k)=0$}
{Set $d_{p_v}(k) \leftarrow d_{p_v}(k) +1 $ \;}
Insert the pair $(d_v,v)$ into the \texttt{PQ}\;
}
}
}
}
\For{each \texttt{reached} $i\in V^\prime$}
{
Set $C^\omega\leftarrow C^\omega +1$\;
\For{each label $k \in K_{p_i}$}
{ 
Set $c_k^\omega\leftarrow c_k^\omega + d_{p_i}(k)$ 
%for $k:(u,v)\in A_k^\omega$
\;}
}
\end{algorithm}

Another characteristic of the implementation with option $\texttt{lift=H}$ is that we carry \MS{out} a pre-processing step, where we enumerate all the seed-to-node pairs having a directed path, which we refer to as \emph{pure label path}, consisting of arcs with the same label. If there is a pure label path from some $i\in I$ to $j\in V^\prime$ with label $k$, then we add an extra arc $(i,j)$ to $A_k$ if it is not already present. This approach supports the lifting procedure and facilitates finding better paths with respect to the objective function 
\eqref{eq:nDSPr_obj_lifted}
although it does not ensure finding the optimal one.

\subsection{Scenario-dependent extended seed sets}

%As described earlier, some of the labels cannot be blocked and they are aggregated under label 0. 
%\footnote{\KT{I think we removed the unblockable label from the problem definition, so we need to introduce it here?}}
While implementing the \BC algorithm, we explicitly consider the \MS{potential} presence of some labels with an extremely high blocking cost. These labels are not feasible to be blocked and therefore are treated as \emph{unblockable} labels.
Under any feasible blocking decision $\bar{x}$, the seed nodes can infect new nodes via arcs with an unblockable label. Therefore, the nodes that are reachable from $I$ on $G^\omega$ via paths consisting of arcs with such labels behave exactly in the same way as the seed nodes in scenario $\omega$. Let $I^\omega$ denote the set of seed nodes extended with these \emph{seed-like} nodes of scenario $\omega$.  
\begin{proposition}
The objective function of \myref{BEN} can be reformulated as
\begin{equation*}
    \frac{1}{|\Omega|} \sum_{\omega \in \Omega} |I^\omega|+\theta^\omega,
\end{equation*}
%if $I$ is replaced with $I^\omega$ in each \SPa{} and \DSPa{}.
if the Benders optimality cuts are modified so that $\theta^\omega$ estimates the number of activated nodes that are not in $I^\omega$.
\label{prop:extendedSeed}
\end{proposition}
To change the quantity that each $\theta^\omega$ estimates in the desired direction, we need to rewrite the objective functions in our primal Benders subproblems as $\sum_{i\in V\setminus I^\omega} y_i^\omega$ and modify the dual subproblems and the optimality cuts accordingly.
This modification results in a possibly a smaller cut constant,
%$\sum_{i\in I^\omega}  \hat{\alpha}_{i}^{\omega}$,
which in turn could lead to tighter cuts due to down-lifting the cut coefficients 
based on the cut constant. Note that $I^\omega$ is independent of $\bar{x}$, and therefore it needs to be computed only once in a pre-processing step.

\subsection{Cut sampling}
\label{section:cutsampling}

In our basic algorithmic setting, we generate one Benders optimality cut for each scenario $\omega \in \Omega$. However this is not mandatory for correctness of the algorithm as one violated cut is enough to cut off an infeasible solution \MS{$(\bar \theta,\bar x)$}. To take advantage of this flexibility, we include the option to sample a subset of $\Omega$ to generate cuts, i.e., sample cuts from the set of all possible cuts.
%Since fractional cuts are optional in our algorithm, whenever we choose to include them in the algorithm we only generate a subset of the cuts. 
We denote by $\tau$ the ratio of the \MS{number of} scenarios \MS{for which} we aim to generate a violated cut for a given solution $(\bar{\theta},\bar{x})$. 
%We first shuffle the scenarios and solve the subproblem for each scenario until we have $\tau_f \times|\Omega|$ violated cuts or all scenarios are considered.
We sort the scenarios in non-decreasing order of the $\bar{\theta}^\omega$ values (ties are broken arbitrarily) and following this order we generate cuts until we have $\tau \times|\Omega|$ violated cuts or all scenarios \NA{have already been} considered. \MS{The cut sampling procedure is detailed in Algorithm \ref{alg:Separation}.}

\begin{algorithm}[htbp]
\caption{\texttt{CutSampling($\bar{\theta}, \bar{x}$,\texttt{lift},$\tau$)}}\label{alg:Separation}
\SetKwInOut{Input}{Input}\SetKwInOut{Output}{Output}
%\Input{A solution ($\bar{\theta},\bar{x}$) to the current \BMP}
%\Output{A (violated cut $\theta^\omega\geq C^\omega - \sum_{k\in K}c^\omega_k x_k$}
Set the desired number of violated cuts $count_{max}=\tau |\Omega|$ and $count\leftarrow 0$ \;
Define $\Omega^\prime \leftarrow \Omega$ \;
\While{$count < count_{max}$ and $\Omega^\prime \neq \emptyset$ }
{
Choose $\omega^\prime = \argmin_{\omega \in \Omega^\prime} \bar{\theta}^\omega$, set $\Omega^\prime \leftarrow \Omega^\prime \setminus \{\omega^\prime\}$\; 
Obtain a cut $BC \leftarrow \texttt{Separation($\omega^\prime, \bar{\theta}, \bar{x}$,\texttt{lift})}$\;
\If{$BC$ is violated at $(\bar{\theta},\bar{x})$}
{
Add $BC$ to BMP, set $count\leftarrow count +1$
}
}
\end{algorithm}

\subsection{Initial solution and initial cuts}
\label{sec:initialization}

%\KT{TODO: rewrite and explain initial cuts (obtained from the greedy solution, one for each scenario)}

%$|K|\times |\Omega|$ initial cuts. For each $k\in \{1,\ldots, |K|\}$, we generate the blocking strategy $\bar{x}$ with $\bar{x}_k=1$ and $\bar{x}_\ell=0$ for $\ell \neq k$. Then we execute for each $\omega\in \Omega$ the separation procedure that we use for integer separation, to generate $|\Omega|$ cuts for $\bar{x}$

%An initial heuristic solution could be fed into the \BC algorithm so that an initial non-trivial primal (upper) bound is available. 
\MS{We implement a greedy heuristic to generate an initial primal solution for our \BC algorithms.}
%A greedy algorithm can be used to obtain such a primal solution although it would not provide an approximation guarantee due to the fact that the objective function is not supermodular, as shown in Section \ref{sec:deftheory}.
In the greedy algorithm, starting with no blockings, i.e. $\mathbf{x}=0$, at each iteration we determine the label $k^\prime$ whose blocking causes the largest marginal reduction in the final spread and set $x_{k^\prime}=1$. We continue until no more blocking is possible within the budget $B$.

In order to obtain an initial non-trivial dual (lower) bound, it is possible to initialize the Benders master problem with a set of optimality cuts. To this end, we propose to use the greedy solution and generate one cut of type \eqref{eq:BendersCut}, \eqref{eq:BendersCut2}, or \eqref{eq:BendersCut3} (depending on the current algorithm setting) for each $\omega \in \Omega$.  

\section{Computational results \label{sec:results}}

The solution algorithm was implemented in C++ and IBM ILOG CPLEX 12.10 was used to solve the MIP models. The computations were made on a single core of an Intel Xeon E5-2670v2 machine with 2.5 GHz and 3GB of RAM, and all CPLEX settings were left on their default values. 

%We set a timelimit 1200 seconds for our all runs. 

\subsection{Instances}
\label{sec:instances}

Our instances \MS{are based on} two types of \NA{networks}: real social networks and randomly generated \NA{ones}. We obtain the real social network data at \url{https://snap.stanford.edu/data/\#socnets}. Two of them, Twitter and Epinions are online social networks, Enron is a communication network, and HepPh is a collaboration network. For random \NA{networks}, we consider the Barab\'asi-Albert (BA)\citep{barabasi1999emergence} and Erd\H{o}s-R\'enyi (ER)\citep{Erdos2022OnRG} models and use the \texttt{R} package \texttt{igraph} to generate them \citep{csardi2006igraph}.
The type of the underlying graphs (directed/undirected), the number of nodes $n$ and the number of arcs/edges $m$ are provided in Table \ref{tab:instances}, together with the number of parallel arcs, if any exist. \MS{As the \Problem is defined on directed graphs, we consider \NA{below} directed versions of the undirected graphs. More detail on the \NA{conversion} procedure is given below, where label creation is also discussed. We observe that most of the instances in Table \ref{tab:instances} do not have any parallel arcs and the instances with parallel arcs just have very few. We thus generate another set of instances based on the instances in Table \ref{tab:instances} by randomly adding (additional) parallel arcs. However, the computational results for these instances are quite similar to the results obtained when considering the instances in Table \ref{tab:instances}. Thus, the results for the instances with (additional) parallel arcs can be found in \ref{appendix_multiple_arcs} for ease of readability.}

\begin{figure}[h!tb]
    \centering
    \begin{subfigure}[b]{0.3\linewidth}
    \includegraphics[width=\linewidth]{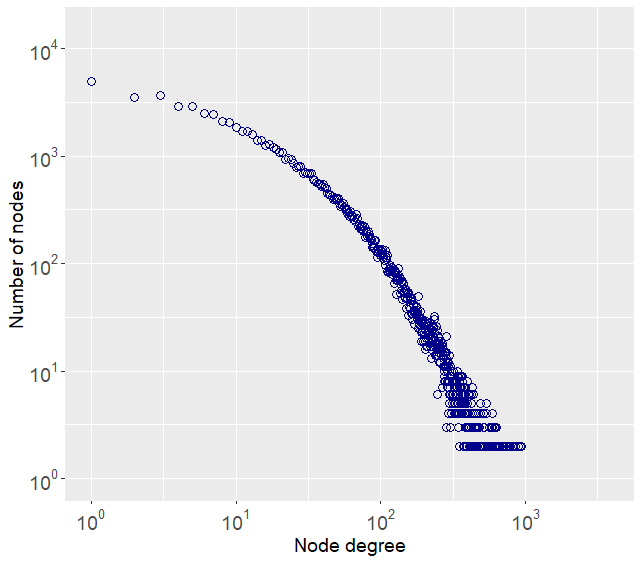}\caption{Twitter}
    \end{subfigure}
    \begin{subfigure}[b]{0.3\linewidth}
    \includegraphics[width=\linewidth]{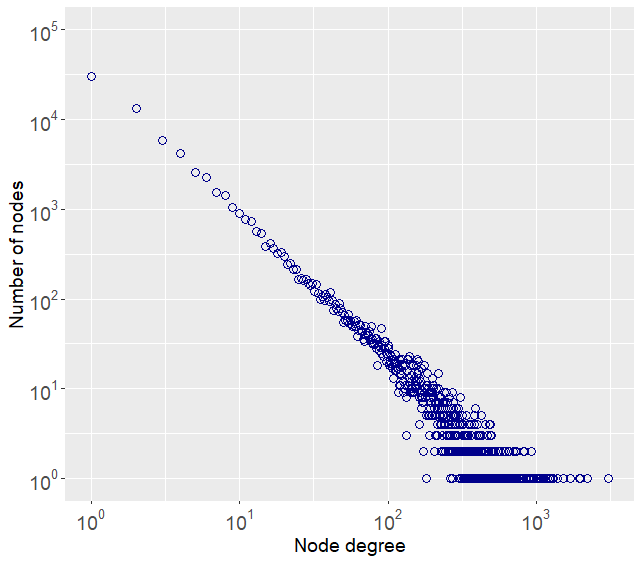}\caption{Epinions}
    \end{subfigure}
    \begin{subfigure}[b]{0.3\linewidth}
    \includegraphics[width=\linewidth]{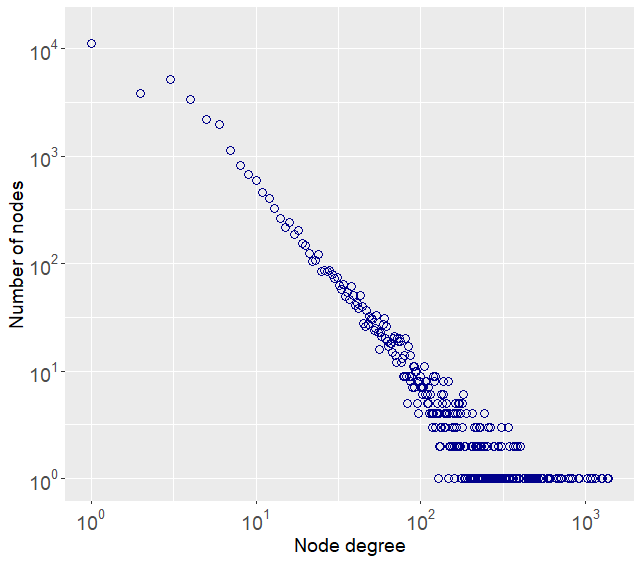}\caption{Enron}%
    \end{subfigure}
    
    \begin{subfigure}[b]{0.3\linewidth}
    \includegraphics[width=\linewidth]{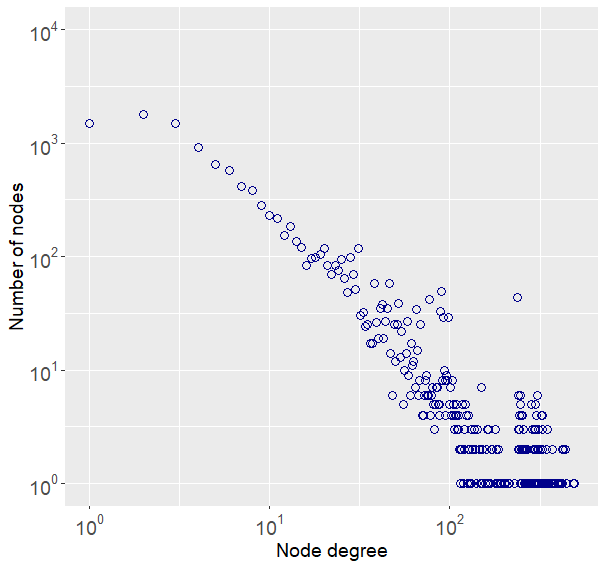}\caption{HepPh}%
    \end{subfigure}
    \begin{subfigure}[b]{0.3\linewidth}
    \includegraphics[width=\linewidth]{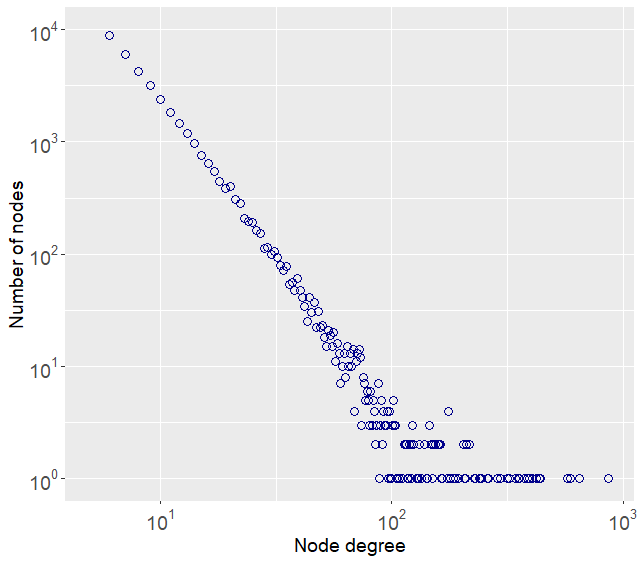}\caption{BA.1}%
    \end{subfigure}
    \begin{subfigure}[b]{0.3\linewidth}
    \includegraphics[width=\linewidth]{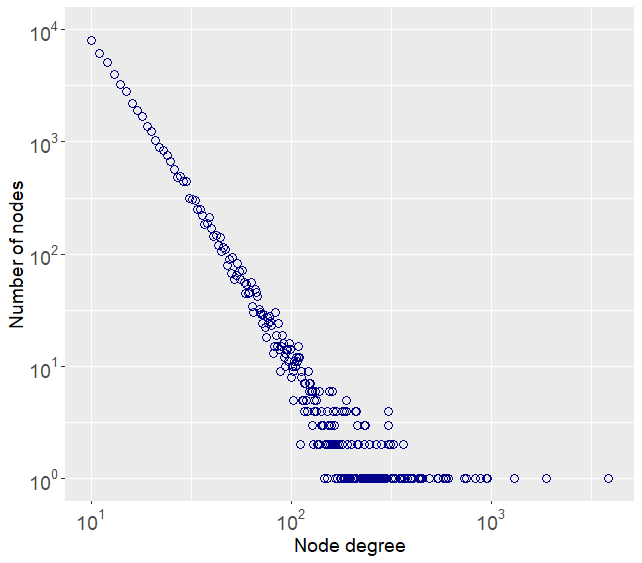} \caption{BA.2}
    \end{subfigure}
    \begin{subfigure}[b]{0.3\linewidth}
    \includegraphics[width=\linewidth]{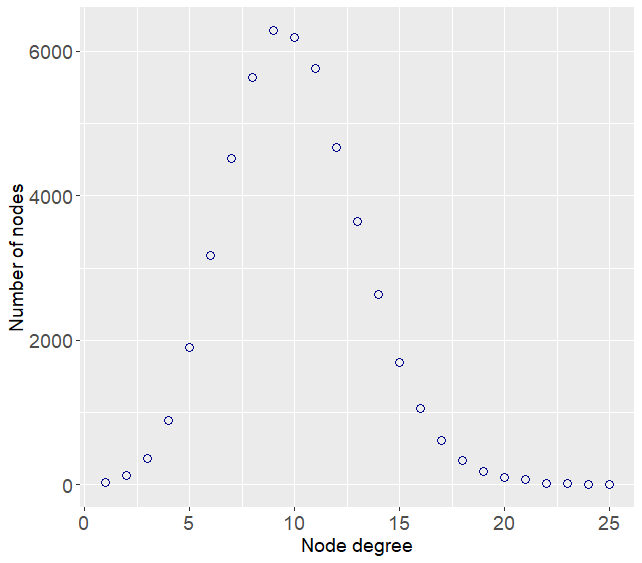} \caption{ER.1}
    \end{subfigure}
    \begin{subfigure}[b]{0.3\linewidth}
    \includegraphics[width=\linewidth]{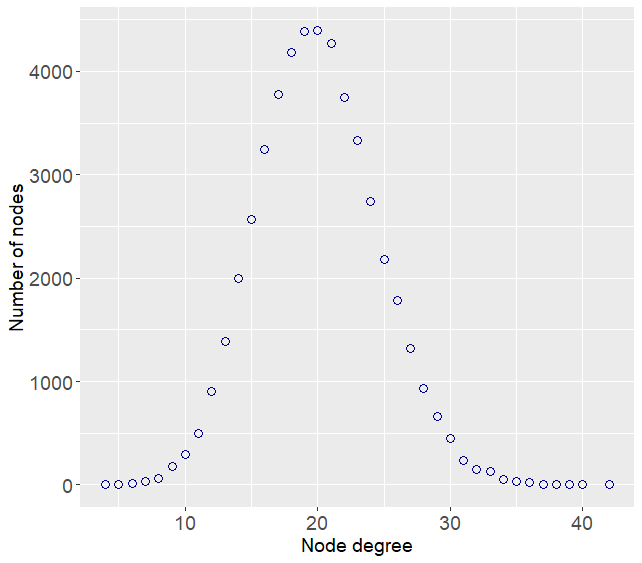} \caption{ER.2}
    \end{subfigure}
    \caption{Degree distributions of the instances}
    \label{fig:degreeDist}
\end{figure}

The node degree distributions of our instances are shown in Figure \ref{fig:degreeDist}. The plots show that, except ER graphs, all the graphs follow a power law distribution.

\begin{table}[htbp]
  \centering \small
  \caption{Description of instances}
    \begin{tabular}{lrrcc}
    \toprule
          & \multicolumn{1}{c}{$n$} & \multicolumn{1}{c}{$m$} & Type & Parallel Arcs\\
    \midrule
    Twitter & 81,306 & 1,768,135 & Directed & $-$\\
    Epinions & 75,879 & 508,837 & Directed & $-$\\
    Enron & 36,692 & 183,831 & Undirected & $-$\\
    HepPh & 12,008 &  118,521 & Undirected & $-$ \\
    BA.1   & 50,000 & 249,995 & Undirected & 310\\
    ER.1   & 50,000 & 250,000 & Undirected & $-$\\
    BA.2    & 50,000 & 499,990 & Undirected & 1807\\
    ER.2    & 50,000 & 500,000 & Undirected & $-$\\
    \bottomrule
    \end{tabular}%
  \label{tab:instances}%
\end{table}%

\iffalse
\begin{table}[htbp]
  \centering \small
  \caption{Description of instances (original and with additional arcs)}
    \begin{tabular}{lcccc}
    \toprule
          & Type & $n$ & $m$ & $m^\prime$ \\
    \midrule
    Twitter & Directed & 81,306 & 1,768,135 & 2,122,399 \\
    Epinions & Directed & 75,879 & 508,837 & 610,720 \\
    Enron & Undirected & 36,692 & 183,831 & 220,719  \\
    HepPh & Undirected & 12,008 &  118,521 &  142,337  \\
    BA.1  & Undirected  & 50,000 & 249,995 & 300,375  \\
    ER.1 & Undirected   & 50,000 & 250,000 & 299,909  \\
    BA.2 & Undirected   & 50,000 & 499,990 & 599,876  \\
    ER.2 & Undirected   & 50,000 & 500,000 & 600,114  \\
    \bottomrule
    \end{tabular}%
  \label{tab:instances}%
\end{table}%
\fi

\paragraph{Assigning arc labels}
We consider $n_L \in \{20,30\}$ labels in addition to one unblockable label $k=0$. We create two classes of instances using the networks described in Table \ref{tab:instances}. For a given network and the value of $n_L$, a label for each arc is generated following a negative binomial distribution with size parameter (the number of targeted successes) equals to one.   
\begin{itemize}
    \item Class 1: The mean of the distribution is chosen as five if $n_L=20$ and eight if $n_L=30$,
    \item Class 2: The mean of the distribution is chosen as eight if $n_L=20$ and twelve if $n_L=30$.
\end{itemize}
Once we sample $m$ numbers following the distribution we choose, we shift the resulting values by +1 and aggregate the ones that are larger than $n_L$ under label 0 so that we have $K = \{0, 1, \ldots, n_L\}$. We assign each number in the resulting vector to the arcs, respectively. 
\KT{Note that as a consequence of this procedure, Class 2 instances have smaller difference between label frequencies.}
If the underlying graph of the original network is undirected, then we assign labels to the edges first, and then replace each edge with two arcs having the same label as the original edge. 

\paragraph{Seed set generation}
We consider \RKT{$|I|\in \{10,50,100\}$} as the size of the seed set $I$. The way we choose $I$ determines the magnitude of the expected spread without any blocking, and the difficulty of limiting the spread. Therefore, we would like to have a seed set for which blocking the spread is not trivial. Since the classical greedy approach \citep{kempe2003maximizing} is inefficient for moderate seed sizes like the ones we consider, we apply a simplified version of the IMM (influence maximization with martingales) algorithm proposed by \citet{tang2015influence}. Differently from the original version of the algorithm, where there is a dedicated sampling phase, we sample a fixed number of 1000 reverse reachable sets (this number is based on preliminary experiments), and then proceed to the node selection phase. Although the algorithm lacks an approximation guarantee with this modification, it still produces good solutions.

\paragraph{Sampling the live-arc scenarios}
For our computational study, we assume that the diffusion follows the IC model. 
Live-arc equivalents of the IC scenarios are obtained by treating each arc $(i,j)$ independently and labelling it as live with a probability equal to the arc probability $p_{ij}$ \citep{kempe2003maximizing}. We end up with a set of live-arcs for each scenario $\omega$ yielding the graph $G^\omega=(V,A^\omega)$. \RKT{We consider three methods while setting up the arc probabilities. In the first and the second one we assume constant probability, $p_{ij}=0.1, \forall (i,j)\in A$ and $p_{ij}=0.05, \forall (i,j)\in A$, respectively, as these are commonly used values}
in related works such as \cite{kempe2003maximizing}, \cite{wu2018two}, \cite{guney2021large}, and \cite{kahr2021benders}. \RKT{In the last one we set $p_{ij}=1/indegree(j)$, which is proposed by \cite{kempe2003maximizing} and assigns a larger weight to a link if its end node does not have many connections.}

In order to determine the number of scenarios to sample, we consider the works on the quality of the approximations via sample average approximation (SAA) and the strategy applied to choose the best sample size. We will briefly mention some of the most significant properties relevant for our case \MS{(see, e.g. \citet{Homemdemello2014} for more details on SAA)}. Let $z^*$ and $z^*_{\Omega}$ be the optimal objective values of the original problem and the approximation problem (\Problem), respectively.
Also, let $x^*$ and $x^*_{\Omega}$ be the optimal solutions (set of labels) provided by the two formulations. Finally, let $(x^*)_\epsilon$ and $(x^*_{\Omega})_\epsilon$ represent the set of all $\epsilon$-optimal solutions to the original and approximation problems, respectively. Since we have a discrete solution set consisting of a subset of all labels and a linear objective function, the following conditions hold for the \Problem:
\begin{enumerate}
	\item The expected value of the optimal objective value of the \Problem is a lower bound on the true optimal objective value, i.e., $\mathbb{E}[z^*_{\Omega}]\leq z^*$
	\item $z^*_{\Omega} \rightarrow z^*$ and $x^*_{\Omega} \rightarrow x^*$ with probability one as $|{\Omega}| \rightarrow \infty$.
	\item The probability of having the optimal SAA solution $x^*_{\Omega}$ as the true optimum solution $x^*$ converges to one exponentially fast as $|{\Omega}| \rightarrow \infty$.
	\item Given the desired significance levels $\rho \in [0,\epsilon)$ and $\alpha \in (0, 1)$, the probability $P \big \{ (x^*_{\Omega})_\epsilon \subset (x^*)_\epsilon \big \} \geq 1-\alpha$. This means that any $\rho$-optimal solution of the \Problem is an $\epsilon$-optimal solution to the original problem with probability at least equal to $1-\alpha$ when the sample size is given as
\begin{align}
|{\Omega}| \geq \frac{3\sigma^2_{\max}}{(\epsilon-\rho)^2}log\left(\frac{\vert \mathcal{X} \vert}{\alpha} \right). \tag{SAA-SS}
\label{samplesize}
\end{align}
\end{enumerate}
In \eqref{samplesize}, $\sigma^2_{\max}$ is a measure of the maximal variance between the set of optimal solutions and the remaining set of solutions. The parameter $\vert \mathcal{X} \vert$ shows the size of the set of all feasible solutions. Thus, the theoretically required sample size $|\Omega|$ may be too large for practical purposes. Therefore, various strategies are provided in the literature for implementing the SAA efficiently. The works
\citet{wu2018two}, \citet{han2018influence},  \citet{taninmics2019influence}, \citet{guney2021large}, \citet{ kahr2021benders}, and \citet{taninmics2022improved} use a single batch of scenarios with a sample size as large as possible. In this work, we follow a similar approach, which is also called the single replication procedure 
%(SRP) 
category of the SAA method \citep{Homemdemello2014}, and consider two values of the number of scenarios, i.e., $|\Omega|\in \{50,100\}$.
\bigskip

In addition to the parameter values explained before, we consider two budget levels $B\in \{4,6\}$ with unit blocking costs, i.e., a cardinality constraint on the blocking strategy. Overall, we obtain \RKT{1152} instances with \RKT{576} instances in each class.

\subsection{Results \label{sec:results2}}
%We \NA{implement} our solution algorithm for the \Problem on our test instances using several settings \NA{which are} summarized as follows.
\MS{We consider the following settings in our computational study.}

\begin{itemize}
    \item LP: The basic implementation of the BD algorithm \MS{based on the arc-based formulation} is carried out, where an LP is solved for each scenario for separation.
    \item I: The combinatorial separation algorithm is used with the option $\texttt{lift}=\texttt{N}$, and only integer solutions are separated.
    \item I+: Within setting I, extended seed sets $I^\omega$ are computed, an initial primal feasible heuristic solution is obtained and initial cuts are added to the Benders master problem, as described in Section \ref{sec:algorithmic_details} 
    \item I+S: In addition to setting I+, \RKT{we apply cut sampling as described in Section \ref{section:cutsampling}}.
    \item I+SF: In addition to setting I+S, we separate fractional solutions \RKT{with cut sampling}.
    \item I+SFP: Differently from setting I+SF we execute the combinatorial separation algorithm with $\texttt{lift}=\texttt{P}$
    \item I+SFH: Differently from setting I+SF we execute the combinatorial separation algorithm with $\texttt{lift}=\texttt{H}$
\end{itemize}

We also consider the greedy algorithm described in Section \ref{sec:initialization} \MS{as a standalone option} and denote it by G. We set the time limit to one hour and run the algorithms until the time limit is reached or a proven optimal solution is obtained.

\paragraph{Preliminary tests}
\RKT{Before presenting the results under all settings mentioned above, we display in Table \ref{tab:tau_val} the results of our preliminary analysis to determine a good value for the cut sampling ratio parameter $\tau$. We use setting I+S which is the most basic setting involving cut sampling, and run experiments on Class 1 instances. Notice that $\tau=1.00$ refers to no sampling, which makes I+S equivalent to I+. The table shows the average solution time in seconds ($t(s)$), the optimality gap percentage (Gap), the number of B\&B nodes explored (nBB), number of optimality cuts added to separate integer solutions (nIntCut), and the number of optimally solved instances out of 576 (nOpt). Due to better time efficiency, we choose $\tau=0.1$ for our computational analysis.

% Table generated by Excel2LaTeX from sheet 'Summary'
\begin{table}[htbp]
  \centering \small
  \caption{Results for different cut sampling ratios}
    \begin{tabular}{rrrrrr}
    \toprule
    \multicolumn{1}{l}{$\tau$} & \multicolumn{1}{l}{$t(s)$} & \multicolumn{1}{l}{Gap} & \multicolumn{1}{l}{nBB} & \multicolumn{1}{l}{nIntCut} & \multicolumn{1}{l}{nOpt} \\
    \hline
    1.00  & 530.7 & 3.0   & 1877.5 & 17778.6 & 539 \\
    0.50  & 395.8 & 1.6   & 2305.1 & 12493.6 & 559 \\
    0.20  & 328.1 & 0.6   & 3015.5 & 8989.5 & 567 \\
    0.10  & 279.2 & 0.4   & 3639.2 & 7889.6 & 568 \\
    0.05  & 285.6 & 0.4   & 3710.4 & 7299.6 & 568 \\
    \hline
    \end{tabular}%
  \label{tab:tau_val}%
\end{table}%

}

\paragraph{Main results}
The aggregated results over all Class 1 (Class 2) instances are shown in Table \ref{tab:Class1} (Table \ref{tab:Class2}) for each of our settings. 
%These tables display the average solution time in seconds ($t(s)$\MS{, where TL indicates that all runs reached the time limit}), the best known upper bound (UB), the best known lower bound (LB), the optimality gap percentage (Gap), lower bound at the B\&B root node ($\text{LB}_0$), number of B\&B nodes explored (nBB), number of optimality cuts added to separate integer solutions (nIntCut), number of optimality cuts added to separate fractional solutions (nFrCut), and number of optimally solved instances out of 128 (nOpt), respectively.
\RKT{In addition to the measures displayed in Table \ref{tab:tau_val}, these tables display the averages of the best known upper bounds (UB), the best known lower bounds (LB), lower bounds at the B\&B root node ($\text{LB}_0$), and the numbers of optimality cuts added to separate fractional solutions (nFrCut), respectively.
}
It is understood from the tables that solving the Benders subproblems in a combinatorial way decreases the solution times dramatically. Moreover, the lower bounds which cannot be improved from a trial value within the time limit under setting LP, becomes much better under setting I. This reflects on the optimality gaps and the number of optimally solved instances. We are able to improve the average gaps further with the other enhancements such as cut sampling and initial cuts. Another big improvement is due to separating fractional solutions which increases nOpt significantly, \RKT{especially for Class 2 instances}. While both lifting options work very well on the overall, 
%the heuristic lifting is not always more efficient than the posterior cut lifting, as the optimality gaps and solution times suggest. The only clear improvement in our measures that the heuristic lifting brings seems to be root lower bounds. The improvement in this measure does not impact the final results except nBB, in the same amount. 
\RKT{the heuristic lifting performs similarly to the posterior cut lifting for Class 1 instances but significantly outperforms it for Class 2 instances. For both instance sets, heuristic lifting yields tighter root bounds and smaller nBB.}
We can see that setting I+SFH manages to solve all instances of Class 1 to optimality, \NA{while} for Class 2, \RKT{14 out of 576} instances remain unsolved with this setting. \KT{The numbers in the tables show that Class 2 instances are more difficult to solve. This can be interpreted as it becomes more difficult to determine the best labels to block when the label frequencies are close.}

% Table generated by Excel2LaTeX from sheet 'Summary'
\begin{table}[htbp]
  \centering \footnotesize
  \caption{Aggregated results for Class 1 instances}
    \begin{tabular}{lrrrrrrrrr}
    \hline
    Method & \multicolumn{1}{l}{$t(s)$} & \multicolumn{1}{l}{UB} & \multicolumn{1}{l}{LB} & \multicolumn{1}{l}{Gap} & \multicolumn{1}{l}{$\text{LB}_0$} & \multicolumn{1}{l}{nBB} & \multicolumn{1}{l}{nIntCut} & \multicolumn{1}{l}{nFrCut} & \multicolumn{1}{l}{nOpt} \\
    \hline
    LP   & 3586.0 & 3678.4 & 134.4 & 83.0  & 54.4  & 209.3 & 5299.6 & 0.0   & 6 \\
    I     & 511.1 & 3339.8 & 3004.5 & 2.8   & 222.9 & 1793.6 & 18857.7 & 0.0   & 542 \\
    I+    & 530.7 & 3320.8 & 2961.8 & 3.0   & 311.2 & 1877.5 & 17778.6 & 0.0   & 539 \\
    I+S   & 279.2 & 3320.8 & 3248.9 & 0.4   & 323.0 & 3639.2 & 7889.6 & 0.0   & 568 \\
    I+SF  & 209.7 & 3320.8 & 3314.3 & 0.0   & 821.3 & 486.2 & 1720.6 & 4649.8 & 575 \\
    I+SFP & 168.3 & 3320.8 & 3320.8 & 0.0   & 1009.8 & 345.8 & 1307.2 & 3271.1 & 576 \\
    I+SFH & 167.2 & 3320.8 & 3320.8 & 0.0   & 1383.5 & 208.2 & 959.2 & 2932.7 & 576 \\
    G     & 29.5  & 3320.9 &       &       &       &       &       &       &  \\
    \hline
    \end{tabular}%
  \label{tab:Class1}%
\end{table}%

% Table generated by Excel2LaTeX from sheet 'Summary'
\begin{table}[htbp]
  \centering \footnotesize
  \caption{Aggregated results for Class 2 instances}
    \begin{tabular}{lrrrrrrrrr}
    \hline
    Method & \multicolumn{1}{l}{$t(s)$} & \multicolumn{1}{l}{UB} & \multicolumn{1}{l}{LB} & \multicolumn{1}{l}{Gap} & \multicolumn{1}{l}{$\text{LB}_0$} & \multicolumn{1}{l}{nBB} & \multicolumn{1}{l}{nIntCut} & \multicolumn{1}{l}{nFrCut} & \multicolumn{1}{l}{nOpt} \\
    \hline
    LP   & 3599.5 & 5368.2 & 149.5 & 87.4  & 61.6  & 206.9 & 5533.4 & 0.0   & 1 \\
    I     & 1123.5 & 5152.4 & 3738.3 & 9.6   & 413.8 & 3013.1 & 34824.9 & 0.0   & 449 \\
    I+    & 1101.9 & 5108.6 & 3763.2 & 9.3   & 712.0 & 3071.0 & 33405.1 & 0.0   & 454 \\
    I+S   & 848.3 & 5108.6 & 4263.2 & 5.2   & 716.1 & 6859.3 & 17944.2 & 0.0   & 487 \\
    I+SF  & 504.9 & 5108.6 & 4664.5 & 2.4   & 1793.3 & 1211.2 & 3055.4 & 8446.4 & 546 \\
    I+SFP & 466.7 & 5109.0 & 4795.5 & 1.6   & 2065.0 & 1031.9 & 2357.4 & 6959.4 & 555 \\
    I+SFH & 400.7 & 5108.6 & 4882.1 & 1.1   & 2414.3 & 736.2 & 1992.3 & 6047.9 & 562 \\
    G     & 36.9  & 5109.7 &       &       &       &       &       &       & \\
    \hline
    \end{tabular}%
  \label{tab:Class2}%
\end{table}%

\begin{figure}[h!tb]
    \centering
    \begin{subfigure}[b]{0.5\textwidth}
         \centering
         \includegraphics[width=\textwidth]{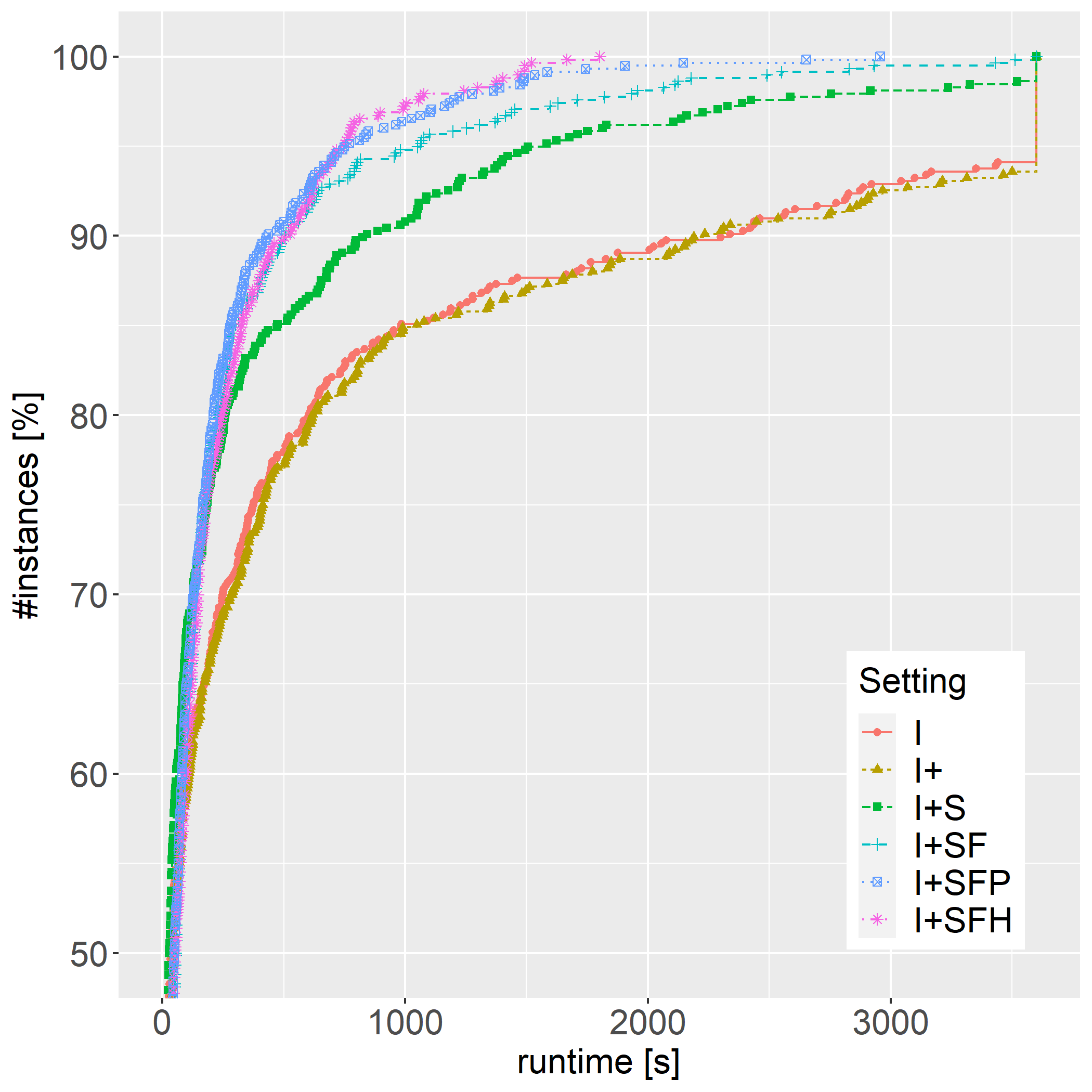}
         \caption{Class 1 running times}
         \label{fig:Class1runtime}
    \end{subfigure}%
    \begin{subfigure}[b]{0.5\textwidth}
         \centering
         \includegraphics[width=\textwidth]{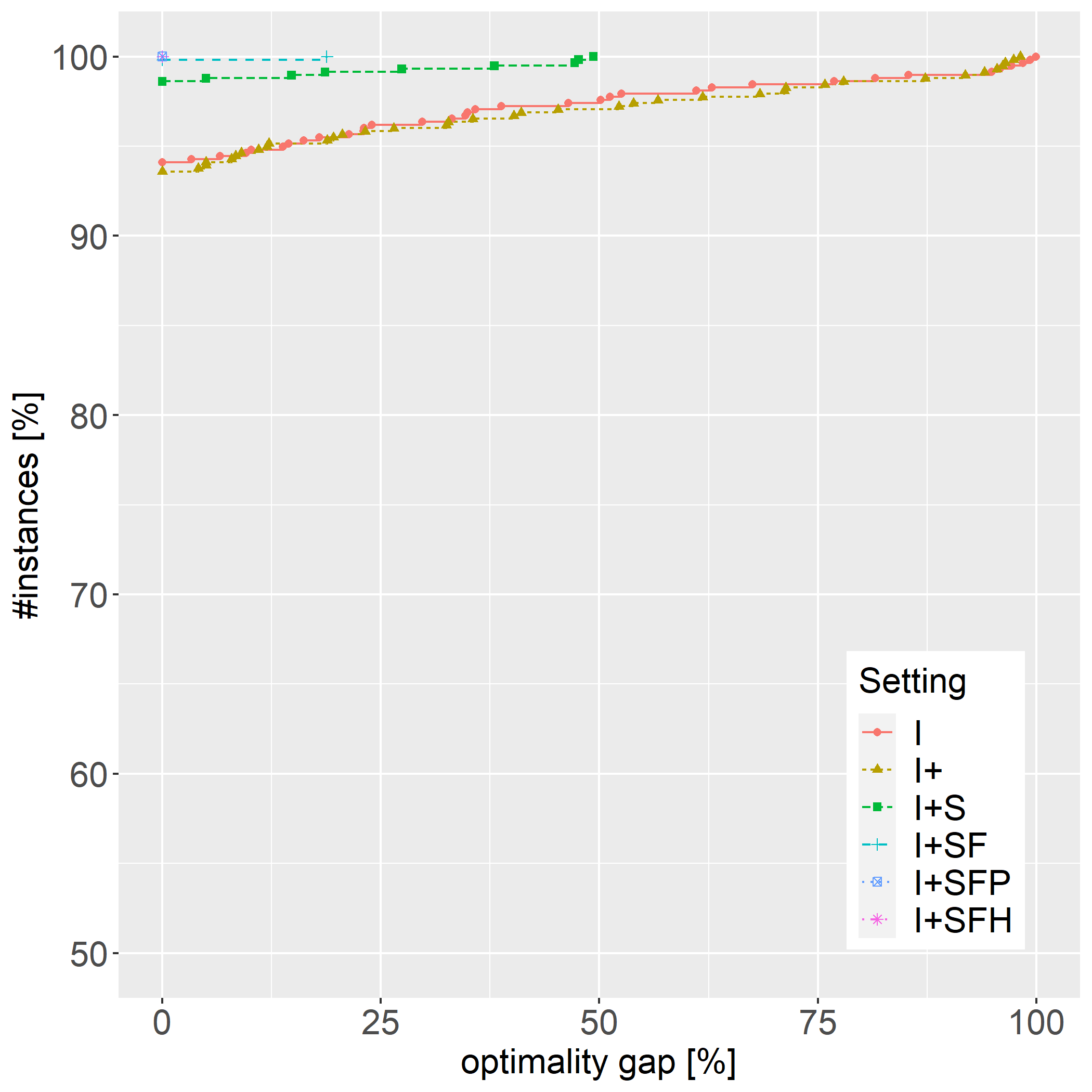}
         \caption{Class 1 final gaps}
         \label{fig:Class2runtime}
    \end{subfigure}
    \begin{subfigure}[b]{0.5\textwidth}
         \centering
         \includegraphics[width=\textwidth]{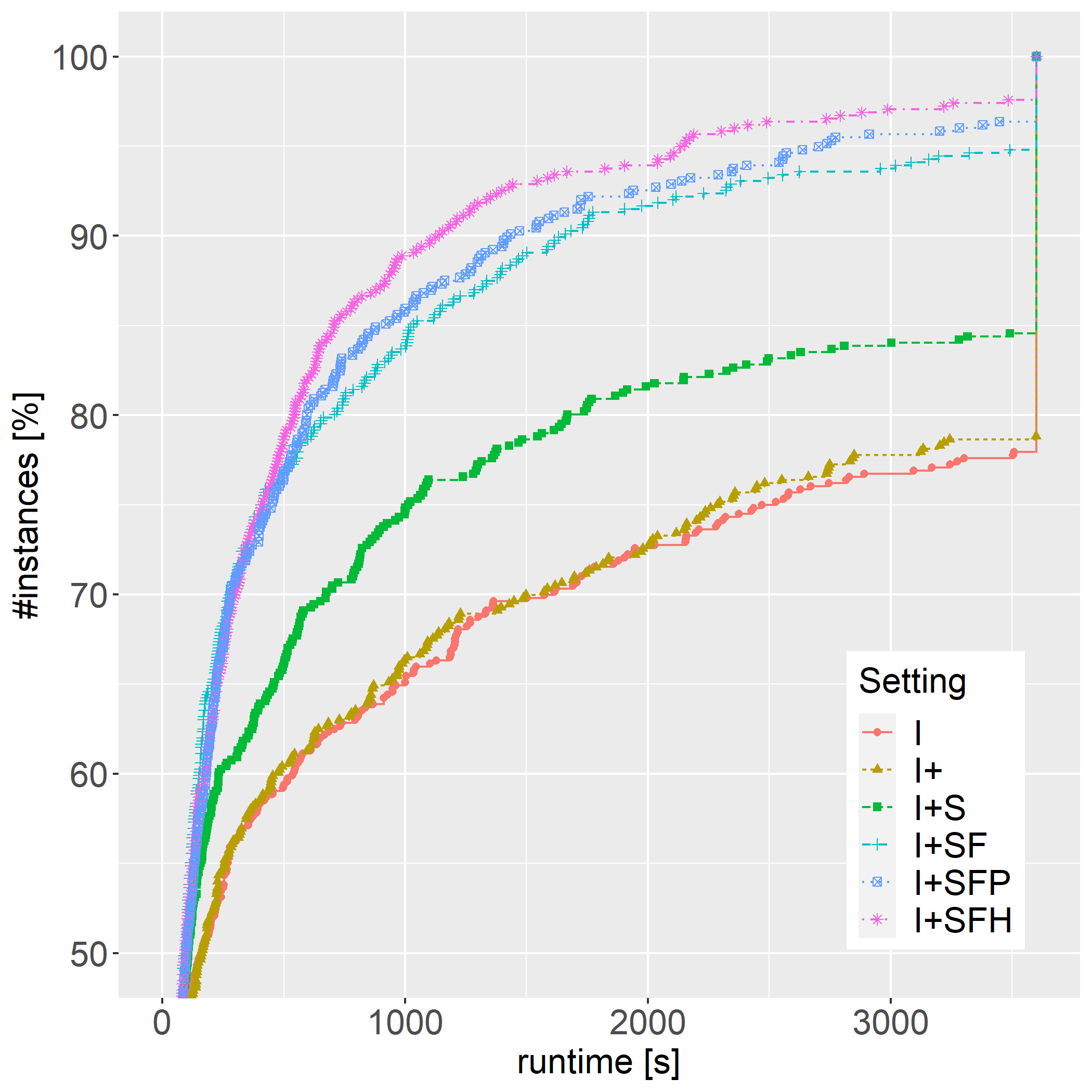}
         \caption{Class 2 running times}
         \label{fig:Class1gap}
    \end{subfigure}%
    \begin{subfigure}[b]{0.5\textwidth}
         \centering
         \includegraphics[width=\textwidth]{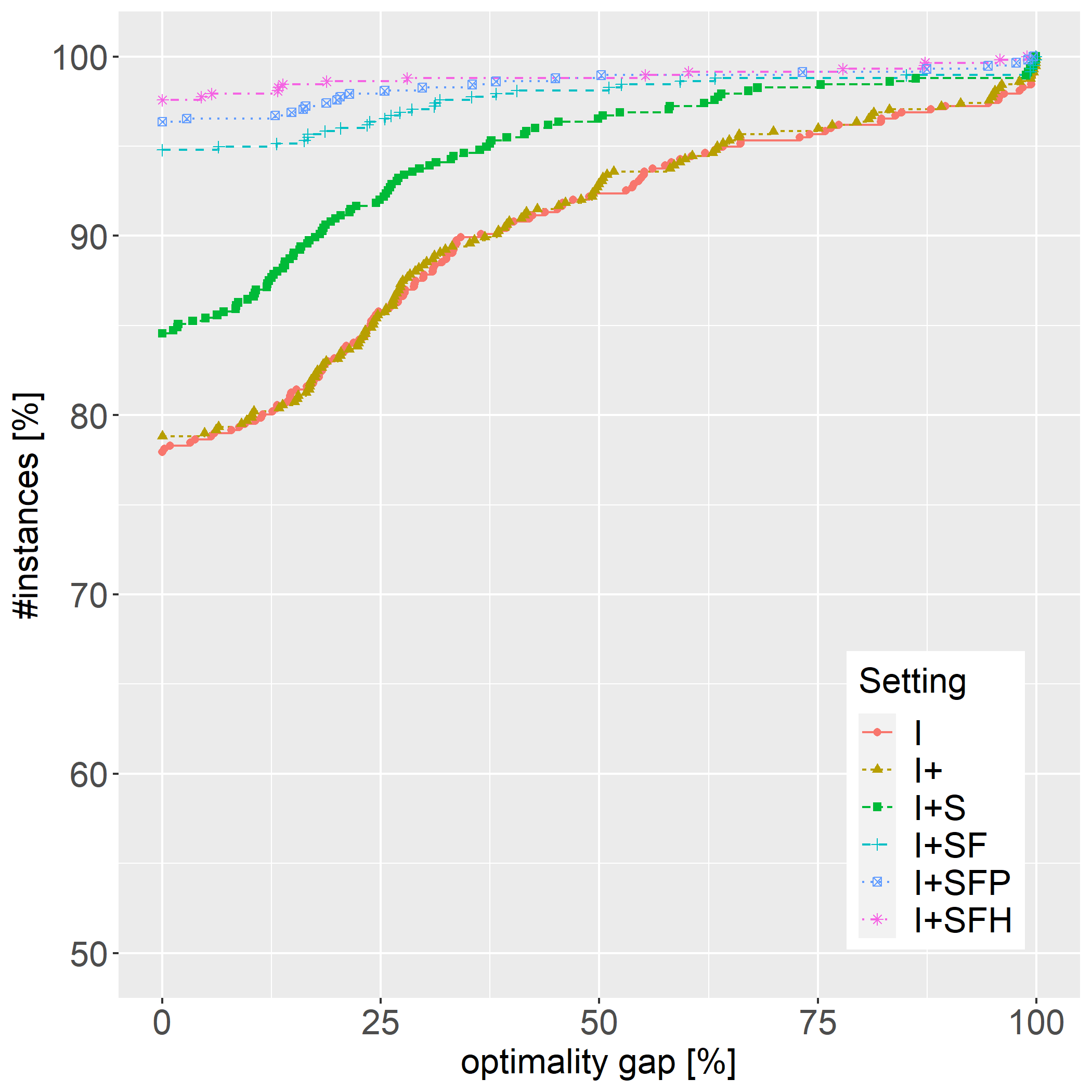}
         \caption{Class 2 final gaps}
         \label{fig:Class2gap}
    \end{subfigure}
    \caption{Cumulative distribution of running times and final gaps for Class 1 and Class 2 instances}
    \label{fig:cum_plots_time}
\end{figure}

\begin{figure}[h!tb]
    \begin{subfigure}[b]{0.5\textwidth}
         \centering
         \includegraphics[width=\textwidth]{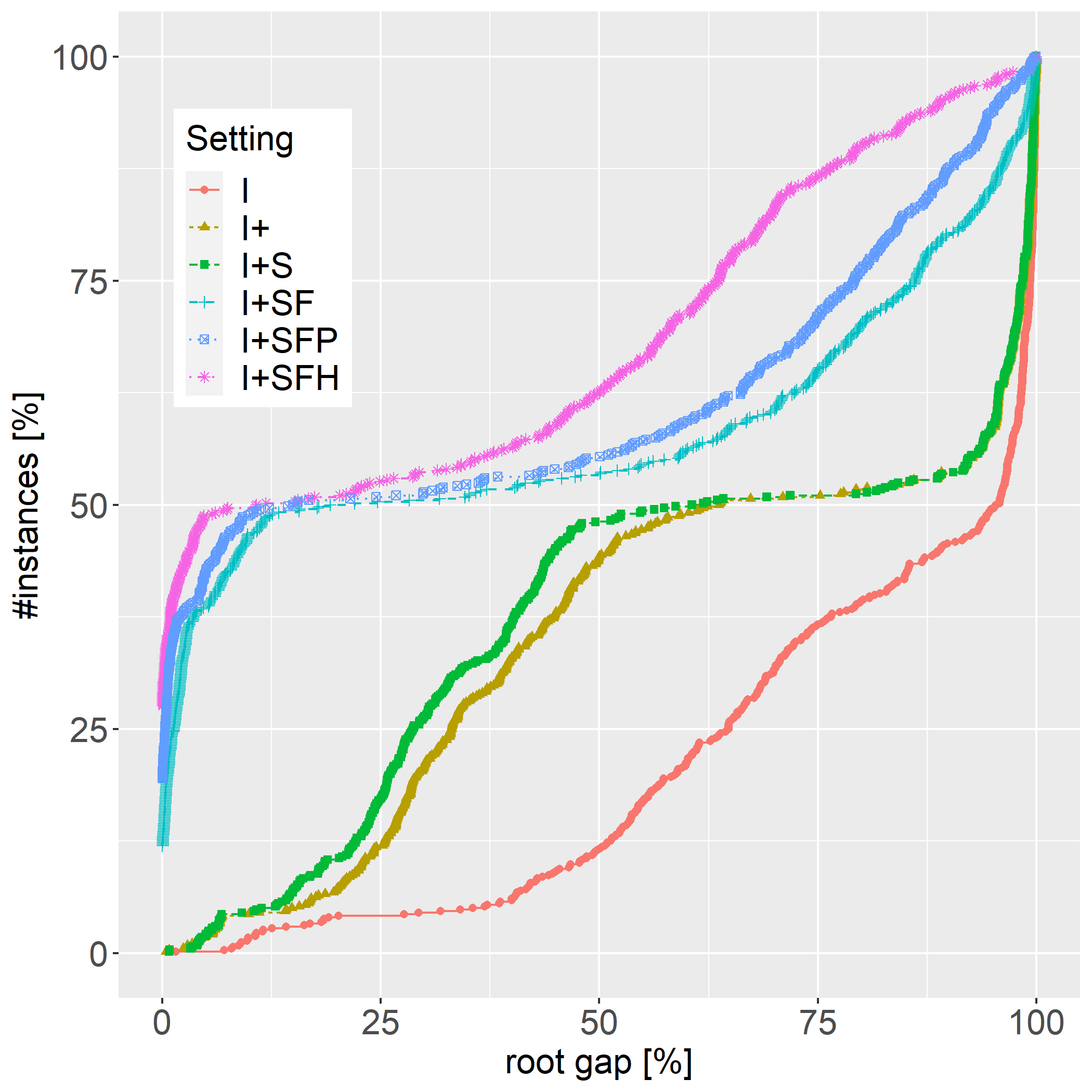}
         \caption{Class 1 root gaps}
         \label{fig:Class1rootgap}
    \end{subfigure}%
    \begin{subfigure}[b]{0.5\textwidth}
         \centering
         \includegraphics[width=\textwidth]{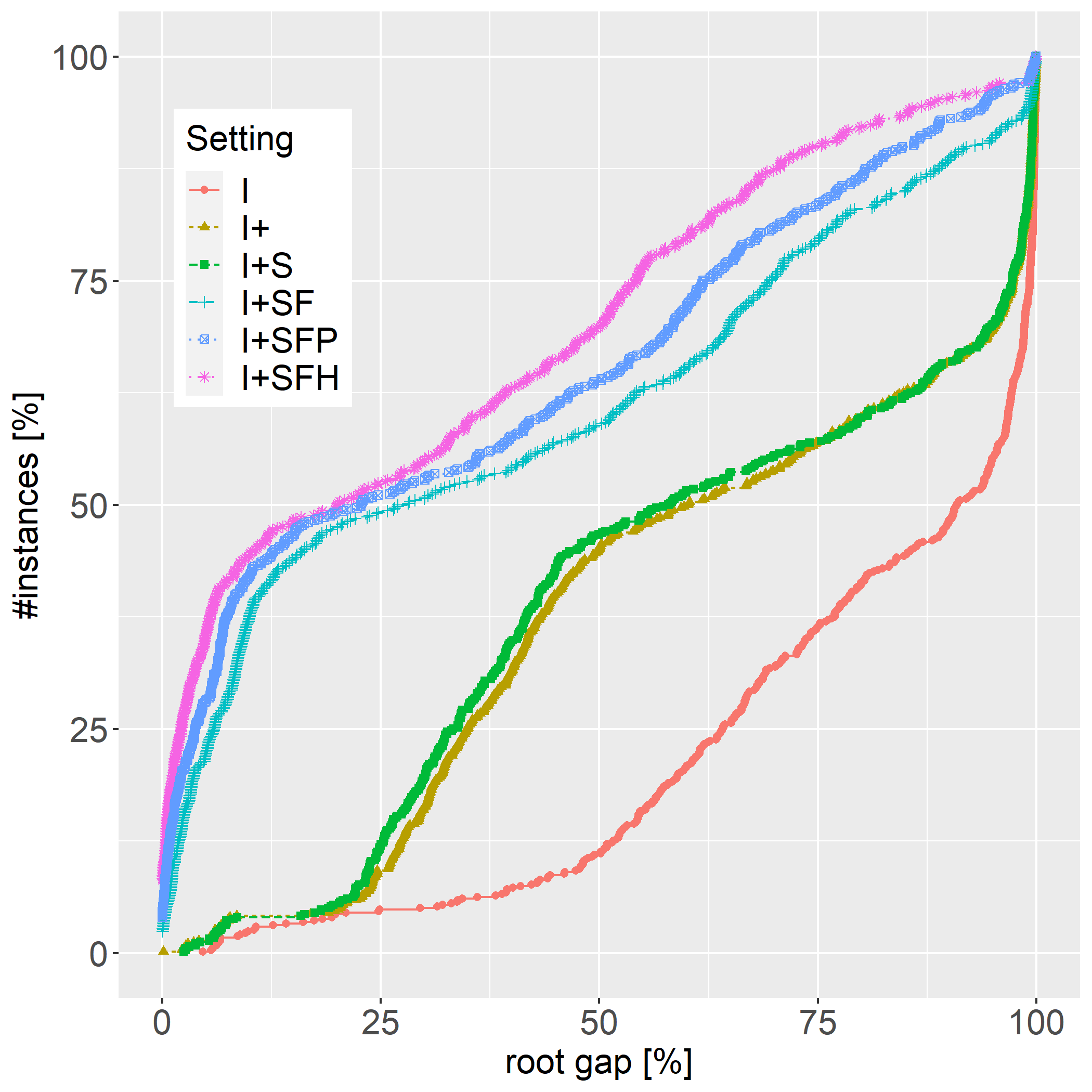}
         \caption{Class 2 root gaps}
         \label{fig:Class2rootgap}
    \end{subfigure}
    \caption{Cumulative distribution root gaps for Class 1 and Class 2 instances}
    \label{fig:cum_plots_root}
\end{figure}

Next, we plot the cumulative distribution \RKT{of the running time and final optimality gaps in Figure \ref{fig:cum_plots_time}, and of the root gaps in Figure \ref{fig:cum_plots_root}}, for all settings except G and LP. The reason is that G is a heuristic and thus does not provide a final gap and takes very short to terminate compared to the exact methods, as expected, and LP takes very long and yields a very large gap. 
%While Figure \ref{} shows the plots for original networks, Figure \ref{} shows the results for the modified networks. 
We see that, the components with largest marginal contribution in terms of runtime are cut sampling (setting I+S) and fractional separation (I+SF), followed by posterior lifting (I+SFP). A similar situation is observed when we focus on the final optimality gaps. In terms of root gaps, initialization (I+), fractional separation, and heuristic lifting (I+SFH) cause the biggest marginal effect.

\paragraph{Impact of the seed set size}

\RKT{Next, we investigate the effect of the seed set size $|I|$ on the performance of our algorithm and its impact on the final spread. We compare the average solution times, objective values, and number of branch-and-bound nodes, as presented in Table \ref{tab:seedSetSize}. The results indicate that the solution times are not significantly affected by increasing the seed set size, nor is the size of the branch-and-bound tree. Moreover, the relative increase in the final spread between $|I|=10$ and $|I|=100$ is around 23\% for instances of Class 1 and 16\% for Class 2.}

% Table generated by Excel2LaTeX from sheet 'Summary'
\begin{table}[htbp]
  \centering \small
  \caption{Average results for low, medium, and large seed set size}
    \begin{tabular}{crrrr}
    \hline
    \multicolumn{1}{l}{Instances} & \multicolumn{1}{l}{$|I|$} & \multicolumn{1}{l}{$t(s)$} & \multicolumn{1}{l}{UB} & \multicolumn{1}{l}{nBB} \\
    \hline
    \multirow{3}[2]{*}{Class 1} & 10    & 160.1 & 2943.5 & 179.5 \\
          & 50    & 177.2 & 3396.3 & 175.0 \\
          & 100   & 184.3 & 3622.7 & 165.0 \\
    \hline
    \multirow{3}[2]{*}{Class 2} & 10    & 413.3 & 4704.2 & 765.5 \\
          & 50    & 389.2 & 5183.2 & 739.1 \\
          & 100   & 399.5 & 5438.4 & 704.0 \\
    \hline
    \end{tabular}%
  \label{tab:seedSetSize}%
\end{table}%

\paragraph{Impact of the arc probabilities}

\RMS{In this paragraph we provide an out-of-sample analysis, i.e., we analyze the quality of an optimal solution (i.e., optimal set of labels to block) obtained for a given scenario $\Omega$ when evaluated under another scenario $\Omega'$.}
\RKT{As explained in Section \ref{sec:instances}, our instances are generated \RMS{by starting with a given network and then using three different arc probability settings, which are used to sample diffusion scenarios $\Omega$ to obtain three different instances}. For a given size of $|I|$, the seed set $I$ of an instance \RMS{are then} determined greedily based on the generated scenarios $\Omega$. Thus, when the problem parameters other than the arc probability setting is fixed (i.e., blocking budget, label distribution, size of the seed set, number of scenarios), we have \RMS{three instances with different $\Omega$ and $I$ for each underlying network}. For a given underlying \RMS{network} (with all parameters except arc probability fixed) let $\Omega_\ell$ denote the scenario set obtained for arc probability setting $\ell$ and $I_\ell$ denote the resulting seed set, for $\ell \in \{1,2,3\}$. To do the out-of-sample analysis, once we solve \Problem for a problem instance with $\Omega_\ell$ and $I_\ell$ for a given $\ell \in \{1,2,3\}$, we obtain the best blocking decision $x^*_\ell$ and then calculate the spread $z^\ell_{\ell'}$ (using $\Omega_{\ell'}$, $I_{\ell'}$, and $x^*_\ell$) for each $\ell' \in \{1,2,3\} \setminus \{\ell\}$.
%Since we use the same seed sets that are used while solving the problem with $\Omega_2$ and $\Omega_3$, we have comparable final spread values. 
Using these values, we calculate $\Delta_\ell = \frac{1}{2}\sum_{\ell' \in \{1,2,3\} \setminus \{\ell\}} 100\times(z^\ell_{\ell'}-z^*_{\ell'})/z^*_{\ell'}$ where $z^*_{\ell'}$ is the best known objective value when the problem is solved with $\Omega_{\ell'}$ and $I_{\ell'}$. This value is the relative deviation of the spread obtained when using the best solution obtained for $\ell$ with the two other arc probability settings, from the best obtained spreads for these settings.
A small value of $\Delta_\ell$ indicates that the solution found under setting $\ell$ performs well for the other two probability settings as well.

In Table \ref{tab:out_of_sample} we report the average $\Delta_\ell$ values over all Class 2 instances that use the probability setting $\ell$, for each $\ell \in \{1,2,3\}$ and for each underlying network. We choose this instance set as the problems in this set are more challenging.
We see that the solutions are usually robust to changes in the arc probabilities, i.e., in $\Omega$. Especially when the problem is solved for $\ell=1$, i.e., $p_{ij}=0.1$ for all $(i,j)\in A$, the resulting solution is either optimal or near-optimal for $\ell \in \{2,3\}$.
The two extreme cases are related to BA.2 and ER.2 instances. In the former, the solutions remain optimal when $\Omega$ change, while in the latter solutions for $\ell=3$ performs on the average $11\%$ worse than the the best known objective value for the probability settings $\ell=1$ and $\ell=2$.
}

% Table generated by Excel2LaTeX from sheet 'out_of_sample'
\begin{table}[htbp]
  \centering \small
  \caption{Deviation of out-of-sample spread from the best known solution, under different arc probability settings}
    \begin{tabular}{lrrr}
    \hline
    \\[-1em]
    Network & \multicolumn{1}{l}{$\overline{\Delta}_1$}& \multicolumn{1}{l}{$\overline{\Delta}_2$} & \multicolumn{1}{l}{$\overline{\Delta}_3$} \\
    \hline
    Twitter & 0.72  & 0.69  & 0.74 \\
    Epinions & 1.02  & 0.56  & 0.61 \\
    Enron & 0.53  & 0.37  & 0.27 \\
    HepPh & 1.27  & 1.80  & 0.78 \\
    BA.1  & 0.65  & 2.05  & 0.68 \\
    ER.1  & 1.50  & 5.44  & 2.73 \\
    BA.2  & 0.00  & 0.00  & 0.00 \\
    ER.2  & 1.12  & 7.29  & 10.98\\
    \hline
    \end{tabular}%
  \label{tab:out_of_sample}%
\end{table}%

\section{Conclusions and outlook} \label{sec:conclusion}

In this paper, we introduce the \emph{measure-based spread minimization problem} (\Problem). The problem is defined on a network with arc labels. It involves selecting a set of arcs labels subject to a budget constraint. For a given selection of labels, all the arcs with these labels are removed from the network. The goal is to select the labels in such a way as to minimize the spread of a contagion, which starts from a given set of seed nodes. We adopt the independent cascade model which is a stochastic diffusion model for the spread dynamics. The \Problem extends the existing studies in the literature that focus on the challenging problem of minimizing the spread of infections by link removal and node deletions by considering a stochastic diffusion model instead of a deterministic one, and by considering removal of labels and \NA{rather than} single node/arc deletions. Moreover, as opposed to many studies, exact solution algorithms are developed to provide a performance guarantee for quite large networks.

We present two integer linear programming formulations for the problem and based on them, we develop Branch-and-Benders-cut solution algorithms. We show that the Benders optimality cuts can be separated in a combinatorial fashion, i.e., without the need of solving linear programs. We also present a set of valid inequalities for one of the formulations and show how these inequalities can be incorporated in the Benders-based algorithms. Moreover, we develop additional enhancements for our algorithms, namely using scenario-dependent extended seed sets, generating initial cuts, and implementing a starting heuristic.
 
There can be several future research directions. First, the model could be extended to implement measures on nodes types in addition to contact types represented by arc labels. For example, all the nodes belonging to a population group could be blocked, which will remove all the associated nodes from network. An example would be to impose a lockdown for people within a certain age interval. Another research direction could be investigating the strategy of partial blocking of arc labels rather than the case of complete blocking in this study. Partial blocking corresponds to allowing people continuing with the contact type but with additional restrictions. For example, an authority may limit the time spent in a restaurant rather than closing it. On the algorithmic side, it could be interesting to try to determine additional valid inequalities. Moreover, a decomposition approach based on Lagrangian decomposition instead of Benders decomposition could also be a fruitful avenue for further work.  

\section*{Acknowledgments}
This research was funded in whole, or in part, by the Austrian Science Fund (FWF)[P 35160-N]. For the purpose of open access, the author has applied a CC BY public copyright licence to any Author Accepted Manuscript version arising from this submission.

%\section*{References}

\bibliographystyle{elsarticle-harv}
%\bibliography{action.bib}

\newpage
\appendix

\section{Proof of Theorem \ref{prop:DSP_solution}}
\label{appendix_proof}

\begin{proof}
We start with showing that $(\hat{\alpha},\hat{\beta})$ is a feasible solution to \myref{DSABF}. By their definitions, $\hat{\alpha}\geq 0$ and $\hat{\beta} \geq 0$. We need to show that $(\hat{\alpha},\hat{\beta})$ satisfies \eqref{eq:DSPr_c1} and \eqref{eq:DSPr_c2}. Consider any $i\in I$. The constraint \eqref{eq:DSPr_c1} can be rewritten as
\begin{equation}
1+\sum_{a\in \delta^{+}_\omega(i)} \beta_a^\omega - \sum_{\delta^{-}_\omega(i)} \beta_a^\omega \geq \alpha_i^\omega .
\label{eq:DSPr_c1repeat}
\end{equation}
At $(\hat{\alpha},\hat{\beta})$, the right-hand-side of \eqref{eq:DSPr_c1repeat} is equal to the number of nodes that are reachable on $G^\omega(\bar{x})$ with an activation path starting at $i$, including itself. By definition of $\hat{\beta}$, $\sum_{a\in \delta^{-}_\omega(i)} \hat{\beta}^\omega_a=0$ since $i$ is a seed node. The value of the first term $\sum_{a\in \delta^{+}_\omega(i)} \beta_a^\omega$ on the left-hand-side is equal to the number of nodes whose activation path contains one of the outgoing arcs of $i$, which is equal to the value of $\alpha_i^\omega$ minus one (note that each seed also activates itself and for this, no arc is needed). Thus, the constraint is satisfied.

Next, consider $i\in V \setminus I$. We show that \eqref{eq:DSPr_c2} is satisfied at $(\hat{\alpha},\hat{\beta})$ by a case distinction. 

\begin{itemize}
\item Case 1: $\sum_{a\in \delta^{-}_\omega(i)} \hat{\beta}^\omega_a = 0$. This means that the total number of nodes that are reachable on $G^\omega(\bar{x})$ via the incoming arcs of $i$, is zero. This implies that $i$ is not reachable as well. Therefore, the number of nodes that are reachable via the outgoing arcs of $i$ is also zero, i.e., $\sum_{a\in A_{i\omega}^{+}} \hat{\beta}^\omega_a=0$. The value of the left-hand-side of \eqref{eq:DSPr_c2} becomes zero and the constraint is satisfied.

\item Case 2: $\sum_{a\in \delta^{-}_\omega(i)} \hat{\beta}^\omega_a > 0$. In this case, at least one of the incoming arcs of $i$ is on the activation path of some node(s) and $i$ is located on some activation path(s). The number of nodes that are reachable via its outgoing arcs, i.e., $\sum_{a\in \delta^{+}_\omega(i)} \hat{\beta}^\omega_a$, is exactly $\sum_{a\in \delta^{-}_\omega(i)} \hat{\beta}^\omega_a-1$ because one of the nodes that are reachable via its incoming arcs is $i$ itself. Thus, the left-hand-side of \eqref{eq:DSPr_c2} is equal to one and the constraint is satisfied. 

\end{itemize}

In order to prove the optimality of $(\hat{\alpha},\hat{\beta})$, we show that the objective value of \myref{DSABF} at $(\hat{\alpha},\hat{\beta})$ is equal to the optimal objective function value of the primal subproblem for the given scenario $\omega$ and solution $\bar x$ of \myref{BEN}. The primal subproblem reads as follows.

\begin{align}
\mytag{PSABF} \quad &\min \sum_{i\in V} y_i^\omega \label{eq:SPr_obj}\\
& \text{s.t.} \notag \\
& \hspace{15pt} y_i^\omega \geq 1 &i&\in I \label{eq:SPr_c1}\\
& \hspace{15pt}y_j^\omega\geq y_i^\omega -\bar{x}_{k} &(&i,j) \in A_k^\omega, k\in K \label{eq:SPr_c2}\\
& \hspace{15pt}y_i^\omega\geq 0 &i&\in V \setminus I
\end{align}

It can be observed that the value of each $y_i^\omega$ for $i \in I$ in an optimal solution to \myref{PABF} is equal to
one due to constraint \eqref{eq:SPr_c1} and it is $\max\{0,1-d_i^\omega(\bar{x})\}$ for $i \in V \setminus I$ due to the cumulative impact of $\bar{x}_k$ in \eqref{eq:SPr_c2}. This results in a total objective value that is equal to the number of $i\in V$ such that $d_i^\omega(\bar{x})<1$, i.e., $i$ is reachable on $G^\omega(\bar{x})$, minus their total shortest path distances (recall that seed nodes $i \in I$ are trivially reachable). 
Now let us consider the objective value of \myref{DSABF} for $(\hat{\alpha},\hat{\beta})$: While $\sum_{i \in I}\hat{\alpha}_i$ accounts for the number of reachable nodes, $\sum_{k\in K} \sum_{a\in A_k^\omega} \hat{\beta}_a^\omega \bar{x}_{k}$ accumulates $\bar{x}_{k}$ for each node whose activation path contains $a$ and as a result keeps track of total shortest path distance. Therefore, we get an objective value that is equal to the optimal objective function value of \myref{PABF}. Together with the feasibility result, this completes the proof.
\end{proof}

\section{Results for instances with (additional) parallel arcs}
\label{appendix_multiple_arcs}

As explained in Remark \ref{remark:multiple}, the underlying graph of the problem instance is allowed to contain parallel arcs between any pair of nodes, where these arcs have distinct labels. In order to observe the effect of these arcs and also the impact of having denser graphs in general, we generate a second set of instances where we created new parallel arcs in addition to the (potentially existing) original ones.
%We also generate instances where there can exist multiple arcs/edges with different labels between the same pair of nodes. 
Towards this end, we take the original instances (i.e., the ones discussed in Section \ref{sec:instances}) and for each existing arc we create one or two copies of the arc with probabilities 0.10 and 0.05, respectively. Each newly created arc is assigned a label randomly such that the original and additional arcs have different labels. The resulting number of arcs are shown under column $m^\prime$ in Table \ref{tab:instances}, next to the original number of arcs $m$.

\begin{table}[htbp]
  \centering \small
  \caption{Description of instances with (additional) parallel arcs}
    \begin{tabular}{lcccc}
    \toprule
          & $m$ & $m^\prime$ \\
    \midrule
    Twitter   & 1,768,135 & 2,122,399 \\
    Epinions   & 508,837 & 610,720 \\
    Enron  & 183,831 & 220,719  \\
    HepPh   &  118,521 &  142,337  \\
    BA.1     & 249,995 & 300,375  \\
    ER.1     & 250,000 & 299,909  \\
    BA.2     & 499,990 & 599,876  \\
    ER.2     & 500,000 & 600,114  \\
    \bottomrule
    \end{tabular}%
  \label{tab:new_instances}%
\end{table}%

\begin{figure}[h!tb]
    \centering
    \begin{subfigure}[b]{0.5\textwidth}
         \centering
         \includegraphics[width=\textwidth]{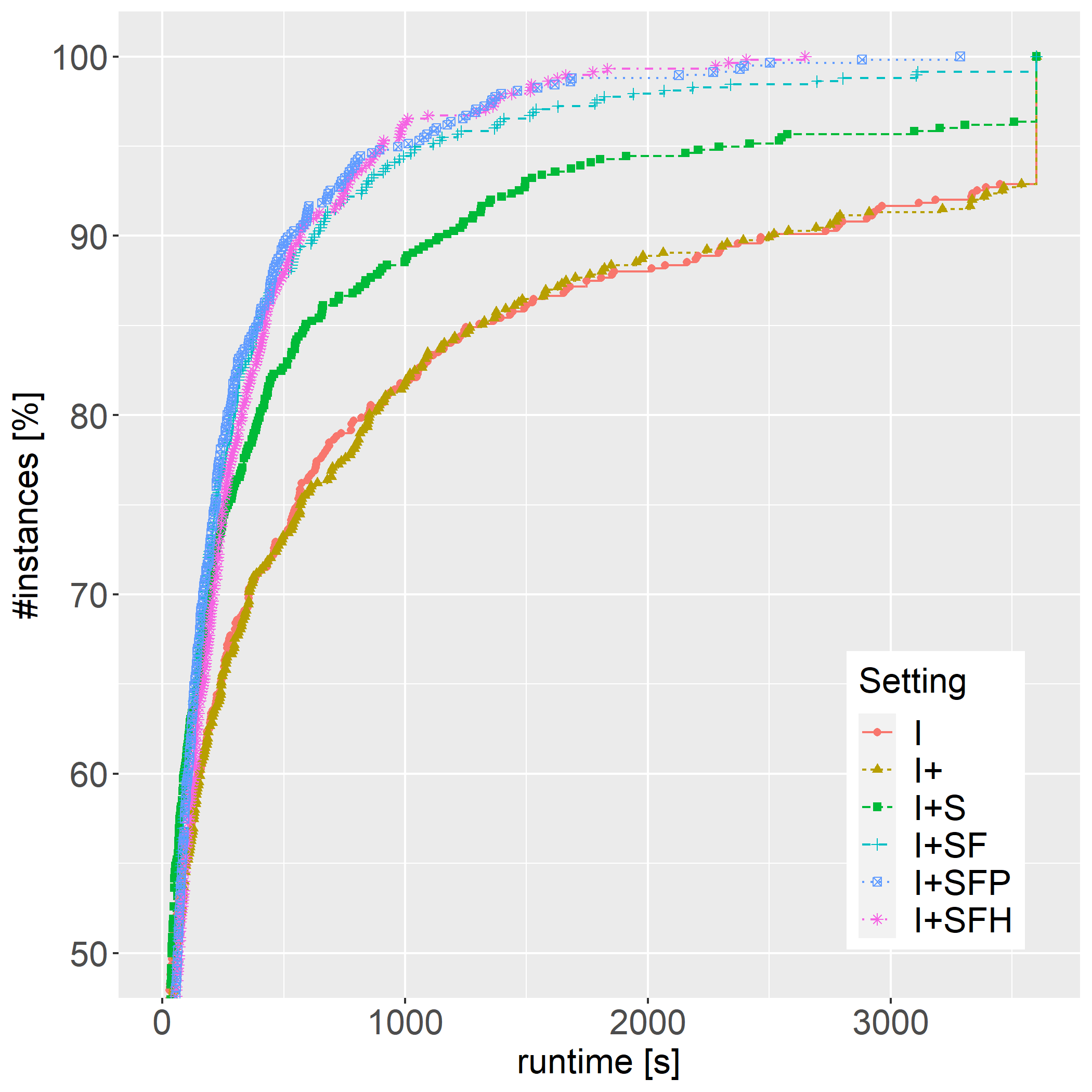}
         \caption{Class 1 running times}
    \end{subfigure}%
    \begin{subfigure}[b]{0.5\textwidth}
         \centering
         \includegraphics[width=\textwidth]{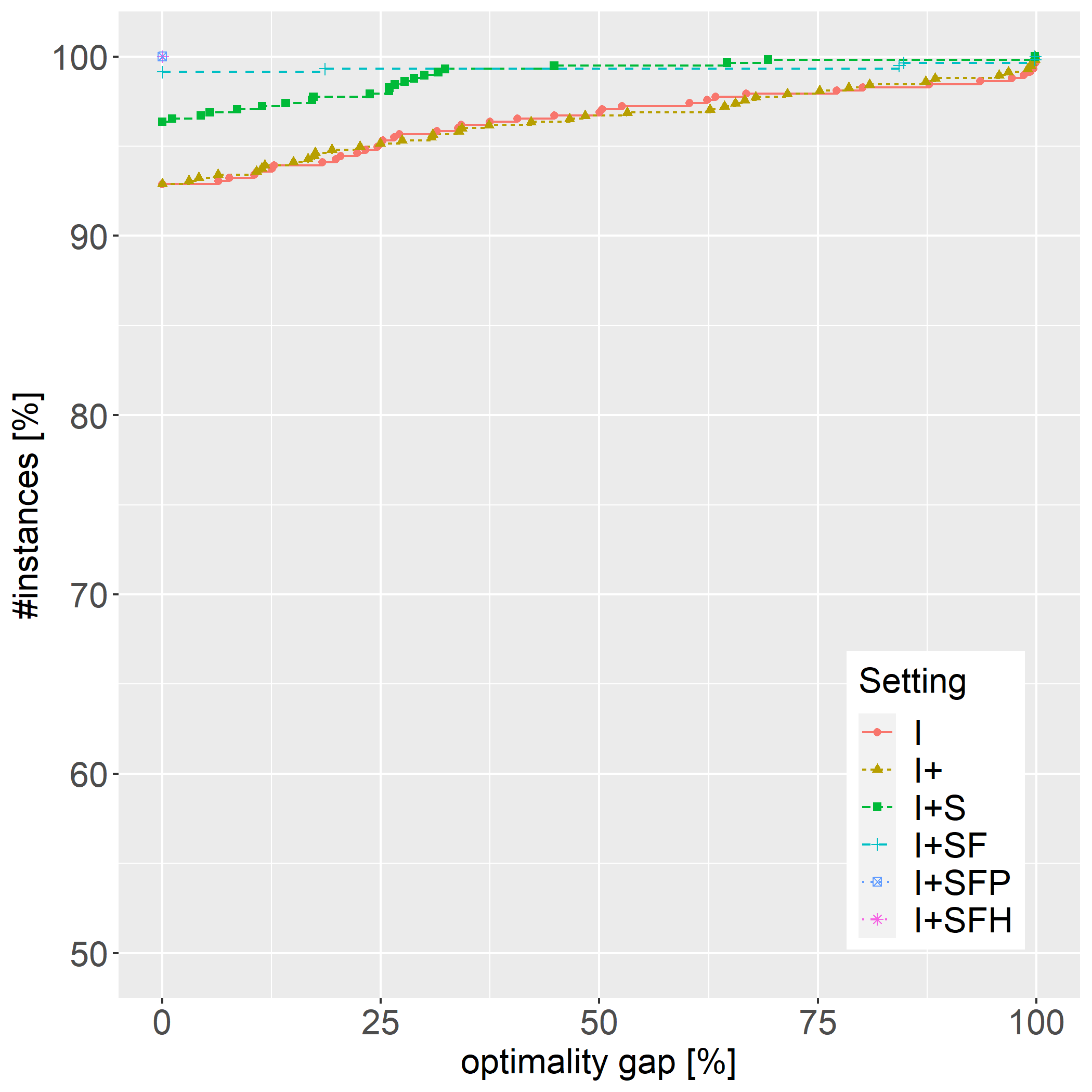}
         \caption{Class 1 final gaps}
    \end{subfigure}
    \begin{subfigure}[b]{0.5\textwidth}
         \centering
         \includegraphics[width=\textwidth]{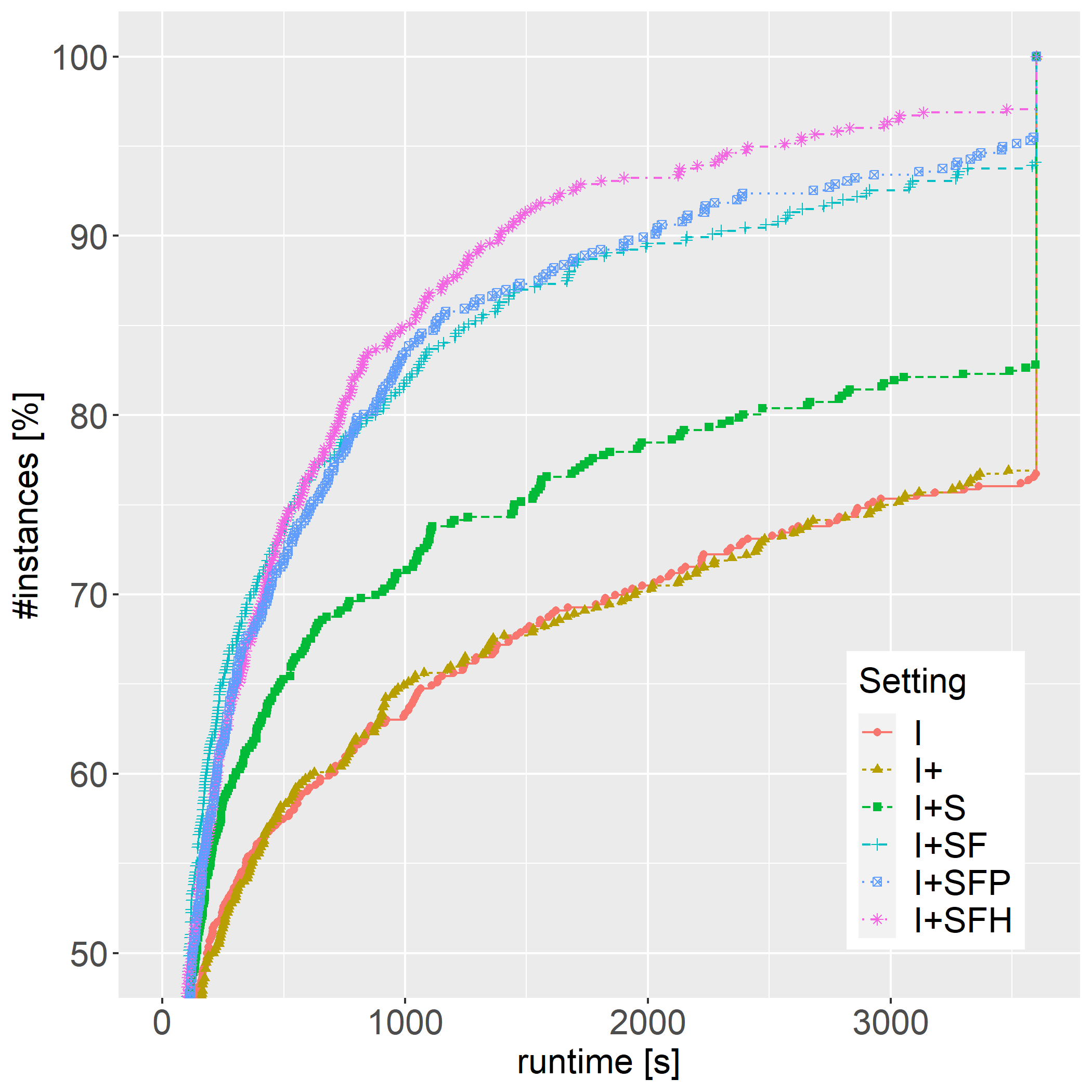}
         \caption{Class 2 running times}
    \end{subfigure}%
    \begin{subfigure}[b]{0.5\textwidth}
         \centering
         \includegraphics[width=\textwidth]{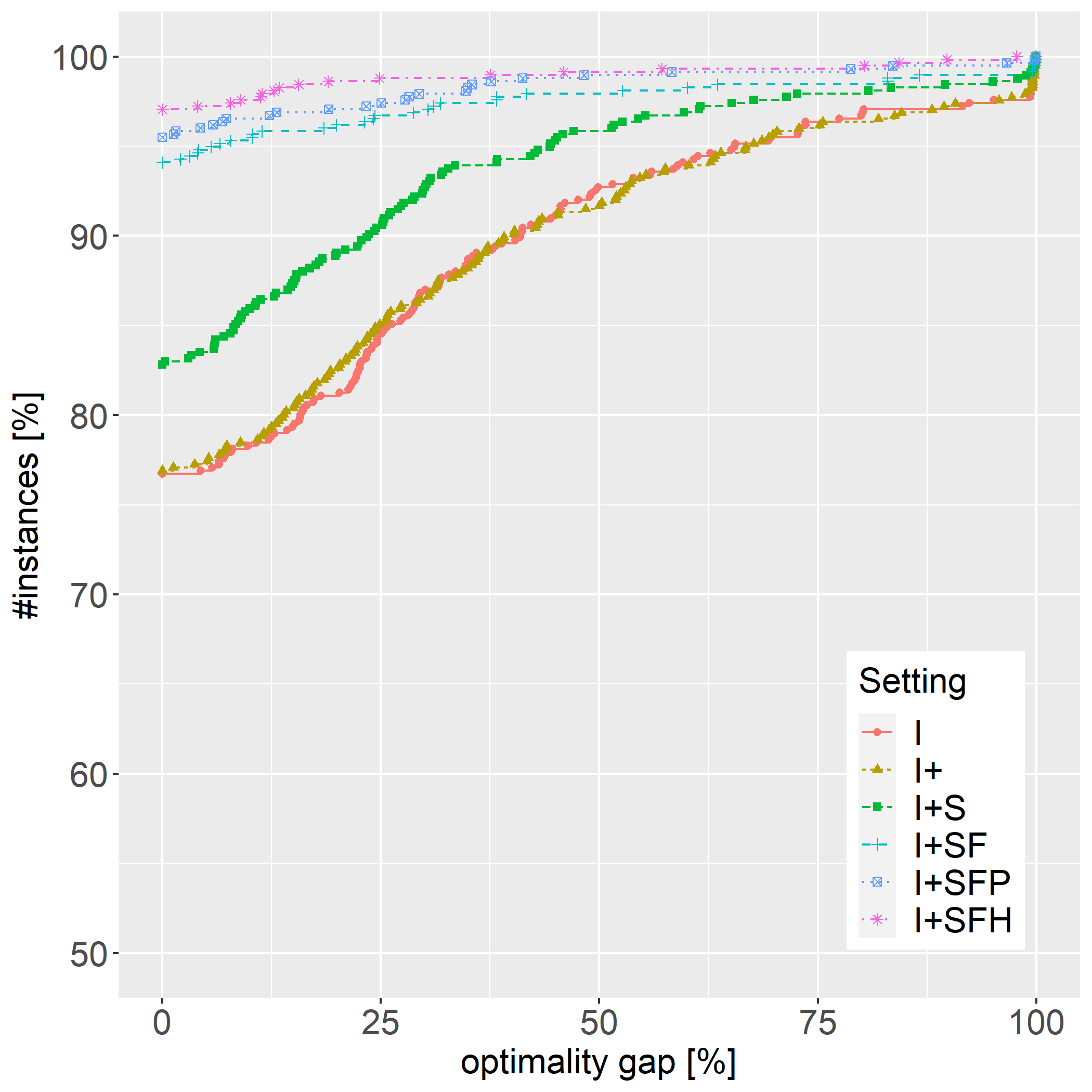}
         \caption{Class 2 final gaps}
    \end{subfigure}
\caption{Cumulative distribution of running times and final gaps for Class 1 and Class 2 instances with (additional) parallel arcs}
    \label{fig:MultArcResults_time}
\end{figure}

\begin{figure}[h!tb]
    \begin{subfigure}[b]{0.5\textwidth}
         \centering
         \includegraphics[width=\textwidth]{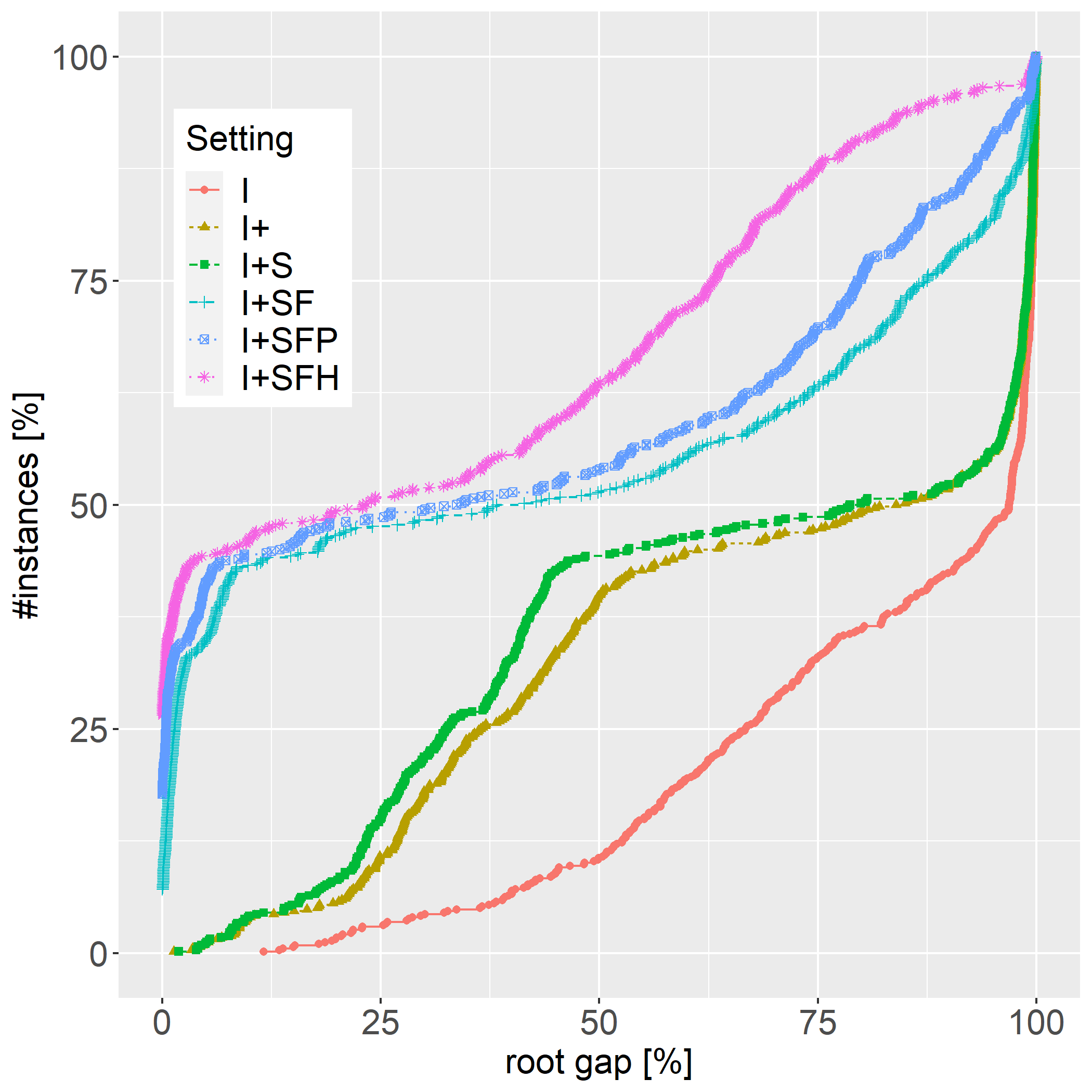}
         \caption{Class 1 root gaps}
    \end{subfigure}%
    \begin{subfigure}[b]{0.5\textwidth}
         \centering
         \includegraphics[width=\textwidth]{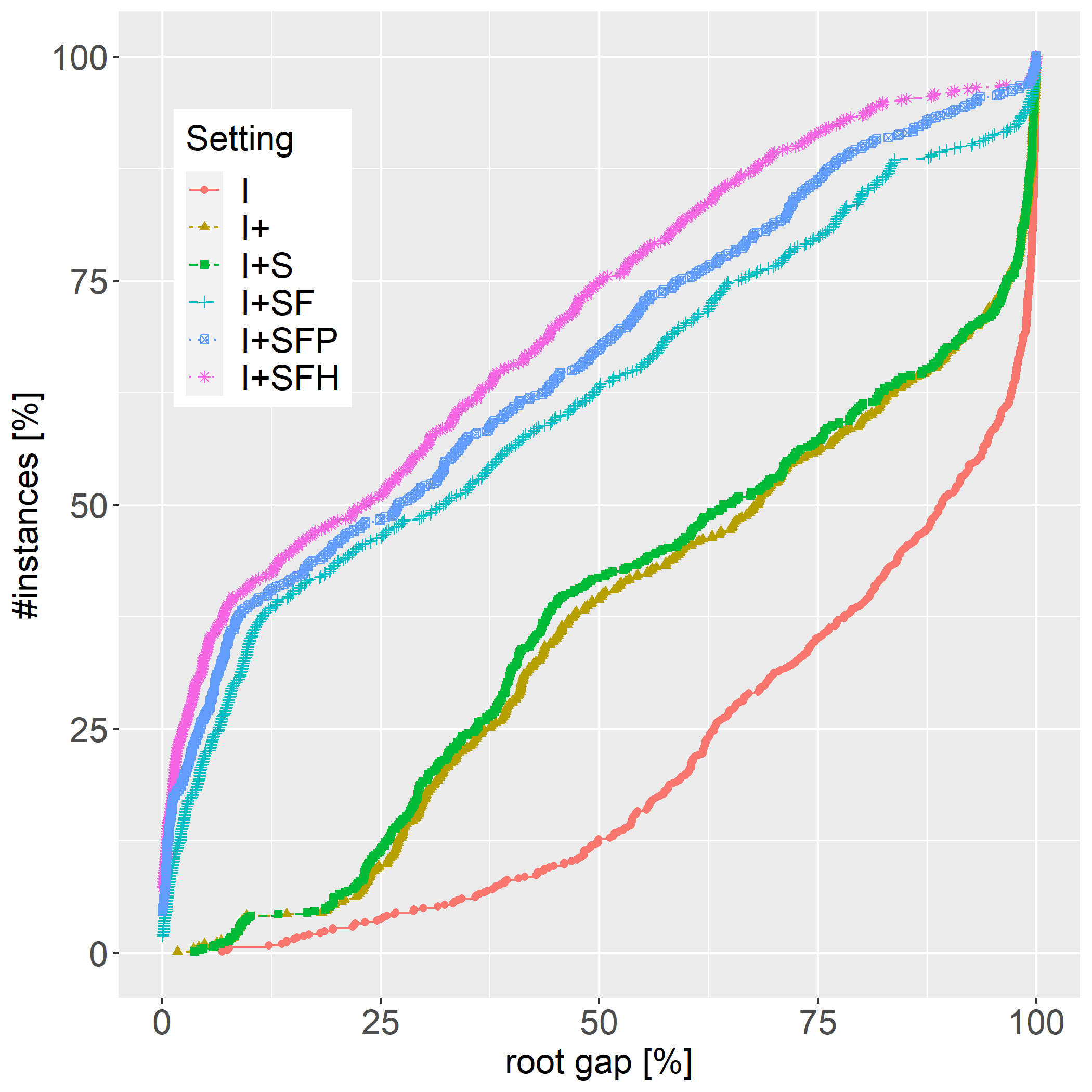}
         \caption{Class 2 root gaps}
    \end{subfigure}
    \caption{Cumulative distribution of root gaps for Class 1 and Class 2 instances with (additional) parallel arcs}
    \label{fig:MultArcResults_root}
\end{figure}

Table \ref{tab:Class1_multiple} (Class 1) and Table \ref{tab:Class2_multiple} (Class 2) present the average results of the algorithms over these new instances. The columns are as explained in Section \ref{sec:results2}. We see that the overall solution times are increased with the newly added arcs/edges. Similarly, the number of optimally solved instances is \RKT{slightly} decreased for both classes. However, the performance of the algorithms with respect to each other remains similar and I+SFP and I+SFH are the best exact algorithms among those we propose, in terms of the running time, final gaps, and nOpt. The cumulative distributions of the running times and final gaps are shown in Figure \ref{fig:MultArcResults_time}, and the cumulative distributions of the
root gaps are shown in Figure \ref{fig:MultArcResults_root}, which are in line with the results of the original instances. Thus, we see that having a significant number of parallel arcs does not seem to have a clear impact on the algorithms.

% Table generated by Excel2LaTeX from sheet 'Summary (paper)'
\begin{table}[htbp]
  \centering \footnotesize
  \caption{Aggregated results for Class 1 instances with (additional) parallel arcs}
    \begin{tabular}{lrrrrrrrrr}
    \hline
    Method & \multicolumn{1}{l}{$t(s)$} & \multicolumn{1}{l}{UB} & \multicolumn{1}{l}{LB} & \multicolumn{1}{l}{Gap} & \multicolumn{1}{l}{$\text{LB}_0$} & \multicolumn{1}{l}{nBB} & \multicolumn{1}{l}{nIntCut} & \multicolumn{1}{l}{nFrCut} & \multicolumn{1}{l}{nOpt} \\
    \hline
    LP   & 3583.9 & 4891.0 & 140.6 & 84.4  & 53.0  & 211.5 & 5174.4 & 0.0   & 7 \\
    I     & 584.4 & 4395.8 & 3736.4 & 3.6   & 218.5 & 1624.1 & 18170.0 & 0.0   & 535 \\
    I+    & 589.5 & 4333.7 & 3711.2 & 3.6   & 297.2 & 1692.2 & 17264.8 & 0.0   & 535 \\
    I+S   & 369.7 & 4333.7 & 4162.3 & 1.1   & 310.1 & 3349.6 & 8049.1 & 0.0   & 555 \\
    I+SF  & 244.1 & 4333.7 & 4264.9 & 0.7   & 912.0 & 455.4 & 1794.7 & 4484.1 & 571 \\
    I+SFP & 207.6 & 4333.7 & 4333.7 & 0.0   & 1131.7 & 337.2 & 1390.7 & 3301.6 & 576 \\
    I+SFH & 214.8 & 4333.7 & 4333.7 & 0.0   & 1716.3 & 193.4 & 982.9 & 2861.4 & 576 \\
    G     & 45.4  & 4333.9 &       &       &       &       &       &       &  \\
    \hline
    \end{tabular}%
  \label{tab:Class1_multiple}%
\end{table}%

% Table generated by Excel2LaTeX from sheet 'Summary (paper)'
\begin{table}[htbp]
  \centering \footnotesize
  \caption{Aggregated results for Class 2 instances with (additional) parallel arcs}
    \begin{tabular}{lrrrrrrrrr}
    \hline
    Method & \multicolumn{1}{l}{$t(s)$} & \multicolumn{1}{l}{UB} & \multicolumn{1}{l}{LB} & \multicolumn{1}{l}{Gap} & \multicolumn{1}{l}{$\text{LB}_0$} & \multicolumn{1}{l}{nBB} & \multicolumn{1}{l}{nIntCut} & \multicolumn{1}{l}{nFrCut} & \multicolumn{1}{l}{nOpt} \\
    \hline
    LP   & 3599.9 & 6845.0 & 146.6 & 88.9  & 66.8  & 189.3 & 5178.5 & 0.0   & 1 \\
    I     & 1190.7 & 6526.7 & 4461.5 & 10.1  & 564.7 & 2694.8 & 34245.1 & 0.0   & 442 \\
    I+    & 1191.5 & 6475.9 & 4467.6 & 10.0  & 918.6 & 2845.2 & 33116.8 & 0.0   & 443 \\
    I+S   & 940.5 & 6475.9 & 5203.9 & 6.1   & 922.4 & 6371.5 & 18890.0 & 0.0   & 477 \\
    I+SF  & 584.4 & 6475.9 & 5853.8 & 2.5   & 2312.3 & 1227.3 & 3207.7 & 8472.5 & 542 \\
    I+SFP & 571.0 & 6475.9 & 6019.7 & 1.7   & 2669.4 & 1061.4 & 2422.9 & 7135.3 & 550 \\
    I+SFH & 481.7 & 6475.9 & 6170.6 & 1.1   & 3096.9 & 760.4 & 2112.6 & 6177.4 & 559 \\
    G     & 51.1  & 6477.9 &       &       &       &       &       &       & \\
    \hline
    \end{tabular}%
  \label{tab:Class2_multiple}%
\end{table}%

% \begin{figure}
%     \centering
%     \includegraphics[width=0.5\linewidth]{New Figures/Mult_Class1_1_final.png}%
%     \includegraphics[width=0.5\linewidth]{New Figures/Mult_Class2_1_final.png}
%     \includegraphics[width=0.5\linewidth]{New Figures/Mult_Class1_2_final.png}%
%     \includegraphics[width=0.5\linewidth]{New Figures/Mult_Class2_2_final.png}
%      \includegraphics[width=0.5\linewidth]{New Figures/Mult_Class1_3_final.png}%
%     \includegraphics[width=0.5\linewidth]{New Figures/Mult_Class2_3_final.png}
%     \caption{Instances with (additional) multiple arcs: Class 1 time, Class 2 time, Class 1 gap, Class 2 gap, Class 1 root gap, Class 2 root gap}
%     \label{fig:MultArcResults}
% \end{figure}

\end{document}